\documentclass[leqno,12pt]{article}
\usepackage{a4wide,amsmath,amssymb,amsfonts,exscale}

\long\def\beginFORGET#1\endFORGET{}

\usepackage[curve,matrix,arrow,cmtip]{xy}

\usepackage{verbatim}

\NoComputerModernTips

\newdir^{ (}{{}*!/-3pt/\dir^{(}}    
\newdir_{ (}{{}*!/-3pt/\dir_{(}}    

\def\Aut{\mathop{\rm Aut}\nolimits}

\def\Frob{\mathop{\rm Frob}\nolimits}
\def\Hom{\mathop{\rm Hom}\nolimits}
\def\Rec{\mathop{\rm Rec}\nolimits} 
\def\CHom{\mathop{\mathcal H\mathit om}\nolimits}

\def\Proj{\mathop{\rm Proj}\nolimits}

\def\Spec{\mathop{\rm Spec}\nolimits}

\def\Quot{\mathop{\rm Quot}\nolimits}

\def\ord{\mathop{\rm ord}\nolimits}

\def\GL{\mathop{\rm GL}\nolimits}
\def\SL{\mathop{\rm SL}\nolimits}

\def\id{{\rm id}}

\newbox\starbox 
\setbox\starbox=\hbox{\vphantom{\underline{\ }}\rlap{\kern4pt\small?}\smash{\lower2pt\hbox{\Large$\square$}}}

\def\hatE{{\mathchoice
  {\hbox{\rlap{\smash{\kern1pt\lower1pt\hbox{$\widehat{\phantom{\hbox{$E$}}}$}}}$E$}}
  {\hbox{\rlap{\smash{\kern1pt\lower1pt\hbox{$\widehat{\phantom{\hbox{$E$}}}$}}}$E$}}
  {\widehat E}
  {\widehat E}}}

\def\hatW{\hbox{\rlap{\smash{\lower1pt\hbox{$\widehat{\phantom{\hbox{$W$}}}$}}}$W$}}

\def\tildeW{\hbox{\rlap{\smash{\lower1pt\hbox{$\widetilde{\phantom{\hbox{$W$}}}$}}}$W$}}

\newbox\checkWbox
\setbox\checkWbox\hbox{\rlap{\smash{\kern.8pt\lower4pt\hbox{\huge \v{}}}}$W$}

\def\circV{{\mathchoice{\circVbig}{\circVbig}{\circVscript}{\circVscriptscript}}}
\def\circVbig{\hbox{\text{\it\r{V}}}}
\def\circVscript{\hbox{\scriptsize\text{\it\r{V}}}}
\def\circVscriptscript{\mbox{\tiny\text{\it\r{V}}}}
\def\circVprime{{\mathchoice{\circV\kern1.8pt{}^\prime}{\circV\kern1.8pt{}^\prime}
                            {\circVscript\kern1.3pt{}^\prime}{\circVscriptscript\kern1pt{}^\prime}}}

\def\circVpprime{{\mathchoice{\circV\kern1.8pt{}^{\prime \prime}}{\circV\kern1.8pt{}^{\prime \prime}}
                            {\circVscript\kern1.3pt{}^{\prime \prime}}{\circVscriptscript\kern1pt{}^{\prime \prime}}}}

\def\bfmu{{\mathchoice{\rlap{$\mu$}\kern.5pt\mu}%
                      {\rlap{$\mu$}\kern.5pt\mu}%
                      {\rlap{$\scriptstyle\mu$}\kern.5pt\mu}%
                      {\rlap{$\scriptscriptstyle\mu$}\kern.5pt\mu}}}                  
\let\phi\varphi
\let\epsilon\varepsilon
\let\setminus\smallsetminus
\let\emptyset\varnothing

\let\le\leqslant
\let\ge\geqslant

\let\geq\geqslant

\newtheorem{Thm}{Theorem}[section]

\newtheorem{Cor}[Thm]{Corollary}
\newtheorem{Def}[Thm]{Definition}

\newtheorem{Exer}[Thm]{Exercise}

\newtheorem{Lem}[Thm]{Lemma}

\newtheorem{Prop}[Thm]{Proposition}

\newtheorem{Rem}[Thm]{Remark}

\renewcommand{\theequation}
             {\arabic{section}.\arabic{Thm}}

\def\qed{{\hskip0pt\unskip\unskip\nobreak\hfil\penalty50
          \hskip1em\hbox{}\nobreak\hfil
          {\bf q.e.d.}%
          \parfillskip=0pt\finalhyphendemerits=0
          \par}\medskip}

\newenvironment{Proof}
               {\noindent{\bf Proof.}\ }
               {\qed}

\newenvironment{Proofof}[1]
               {{\bf Proof of #1:}\ }
               {\qed}

\newcommand{\BA}{{\mathbb{A}}}

\newcommand{\BF}{{\mathbb{F}}}

\newcommand{\BP}{{\mathbb{P}}}
\newcommand{\BQ}{{\mathbb{Q}}}

\newcommand{\BZ}{{\mathbb{Z}}}

\newcommand{\Fa}{{\mathfrak{a}}}

\newcommand{\Fm}{{\mathfrak{m}}}

\newcommand{\CD}{{\cal D}}
\newcommand{\CE}{{\cal E}}
\newcommand{\CF}{{\cal F}}

\newcommand{\CI}{{\cal I}}

\newcommand{\CL}{{\cal L}}

\newcommand{\CO}{{\cal O}}

\def\longto{\longrightarrow}
\def\into{\hookrightarrow}
\let\onto\twoheadrightarrow

\def\longinto{\lhook\joinrel\longrightarrow}

\newcommand{\myatop}{\genfrac{}{}{0pt}{}}

\newbox\mybox
\def\arrover#1{\mathrel{
       \setbox\mybox=\hbox spread 1.4em
              {\hfil$\scriptstyle#1\vphantom{g}$\hfil}
       \vbox{\offinterlineskip\copy\mybox
             \hbox to\wd\mybox{\rightarrowfill}}}}
\def\larrover#1{\mathrel{
       \setbox\mybox=\hbox spread 1.4em{\hfil$\scriptstyle#1$\hfil}
       \vbox{\offinterlineskip\copy\mybox
             \hbox to\wd\mybox{\leftarrowfill}}}}

\def\ontoover#1{\mathrel{
       \setbox\mybox=\hbox spread 1.4em{\hfil$\scriptstyle#1$\hfil}
       \vbox{\offinterlineskip\copy\mybox
             \hbox to\wd\mybox{\rightarrowfill\hskip-2.8mm
                               $\rightarrow$}}}}
\def\leftontoover#1{\mathrel{
       \setbox\mybox=\hbox spread 1.4em{\hfil$\scriptstyle#1$\hfil}
       \vbox{\offinterlineskip\copy\mybox
             \hbox to\wd\mybox{$\leftarrow$\hskip-2.8mm
                               \leftarrowfill}}}}

\newbox\invlimsymbol
\setbox\invlimsymbol=\vtop{\hbox{\rm lim}\vskip-8pt
        \hbox{\hskip1pt$\scriptstyle\longleftarrow$}\vskip-1pt}

\newbox\dirlimsymbol
\setbox\dirlimsymbol=\vtop{\hbox{\rm lim}\vskip-8pt
        \hbox{\hskip1pt$\scriptstyle\longrightarrow$}\vskip-1pt}

\begin{document}


\title{Compactification of a Drinfeld \\ 
Period Domain over a Finite Field}

\author{Richard Pink$^1$, \ \ Simon Schieder$^{2,3}$}
\date{\today}
\maketitle

\footnotetext[1]{Dept. of Mathematics, ETH Z\"urich, 8092 Z\"urich, Switzerland}
\footnotetext[2]{Dept. of Mathematics, Harvard University, Cambridge, MA 02138, USA}
\footnotetext[3]{Supported by the International Fulbright Science and Technology Award of the U.S. Department of State}


\begin{abstract}
We study a certain compactification of the Drinfeld period domain over a finite field which arises naturally in the context of Drinfeld moduli spaces. Its boundary is a disjoint union of period domains of smaller rank, but these are glued together in a way that is dual to how they are glued in the compactification by projective space. This compactification is normal and singular along all boundary strata of codimension~$\ge2$. 
We study its geometry from various angles including the projective coordinate ring with its Hilbert function, the cohomology of twisting sheaves, the dualizing sheaf, and give a modular interpretation for it. We construct a natural desingularization which is smooth projective and whose boundary is a divisor with normal crossings.
We also study its quotients by certain finite groups.
\end{abstract}

\bigskip\bigskip\bigskip
{\bf Mathematics Subject Classification (MSC 2010):} 

Primary: 
14M27	

Secondary:
11F52,	
11T60,	
14G15	

\newpage
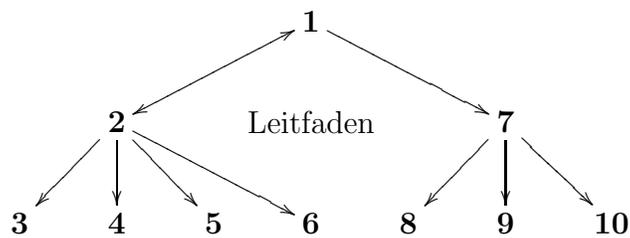
\tableofcontents

\bigskip\bigskip\bigskip\bigskip
$$\xymatrix{
&&& \bf1 \ar[rrd] \ar@{<->}[lld] &&& \\
& \bf2 \ar[dl] \ar[d] \ar[dr] \ar[drr] 
&& \hbox to 0pt{\hss Leitfaden\hss} && \bf7 \ar[dl] \ar[d] \ar[dr] & \\
\bf3 & \bf4 & \bf5 & \bf6 & \bf8 & \bf9 & \bf10 \\}$$

\parindent=0pt
\parskip=\smallskipamount


\newpage

\addtocounter{section}{-1}
\section{Introduction}
\label{Intro}

Let $\BF_q$ be a finite field with $q$ elements. For any positive integer $r$ let $\Omega_r$ be the open dense subscheme of projective space $P_r := \BP^{r-1}_{\BF_q}$ obtained by removing all proper $\BF_q$-rational linear subspaces. This is an interesting algebraic variety over~$\BF_q$ with an action of the finite group $\GL_r(\BF_q)$. By analogy with the (rigid analytic) Drinfeld upper half space associated to a non-archimedean local field in place of $\BF_q$ it has been called a `period domain' by Rapoport \cite{Rapoport-perioddomains} and Orlik \cite{Orlik-Thesis}. It arises naturally as a moduli space of Drinfeld $\BF_q[t]$-modules of rank $r$ with a level structure of level $(t)$. As such, it possesses a natural compactification $Q_r$ analogous to the Satake compactification of Siegel moduli space, which can be characterized using the modular interpretation and/or using Drinfeld modular forms. It turns out that $Q_r$ differs fundamentally from the tautological compactification~$P_r$.

\medskip
The purpose of this paper is to study $Q_r$ as an algebraic variety in its own right from various points of view. We define and analyze it without reference to Drinfeld modules or Drinfeld modular forms---the consequences for these will be explained in the forthcoming paper~\cite{BreuerPink2}. We believe that $Q_r$ carries enough interesting geometry to justify that approach.

\medskip
The basic definitions and some first results are given in Section~\ref{RV1}. We formulate them in a coordinate free manner in order to exhibit the functorial behavior and the action of $\GL_r(\BF_q)$. Thus we will have $\Omega_r = \Omega_V$ and $P_r = P_V$ and $Q_r = Q_V$ for the standard vector space $V=\BF_q^r$. 
We define $Q_V$ by giving a projective coordinate ring $R_V$ for it. In Section \ref{RV2} we present $R_V$ by generators and relations, prove that it is a Cohen-Macaulay normal integral domain, and determine its Hilbert function. These results will be applied in \cite[Sect.$\,$7]{BreuerPink2}, where $R_V$ will be identified with a certain ring of Drinfeld modular forms.

\medskip
In Section \ref{invar} we determine the subring of invariants in $R_V$ under the group $G:=\GL_r(\BF_q)$, under the group $G':=\SL_r(\BF_q)$, and also under a maximal unipotent subgroup $U\subset G$. As a consequence, we show that the quotient varieties of $P_V$ and $Q_V$ under $G$, $G'$, $U$ are all weighted projective spaces of explicitly given weights. Interestingly, although $Q_V$ is in general singular, the quotient $Q_V/U$ is always isomorphic to $\BP^{r-1}_{\BF_q}$ and thus smooth. In Section \ref{URep} we determine the Hilbert function of the ring of invariants in $R_V$ under an arbitrary unipotent subgroup of~$G$.

\medskip
In Section \ref{cohom} we calculate $\dim H^i(Q_V,\CO(n))$ for all integers $i$ and~$n$. In particular we show that it vanishes if $i\not=0$,~$r-1$, as for projective space. In Section \ref{dual} we determine the dualizing sheaf on~$Q_V$.

\medskip
Next we consider the natural stratification of $P_V \cong \BP^{r-1}_{\BF_q}$ whose strata are the $\BF_q$-rational linear subspaces with all smaller $\BF_q$-rational linear subspaces removed. These strata are canonically isomorphic to $\Omega_{V''}$ for all non-zero quotients $V''$ of~$V$. In Section \ref{strat} we show that $Q_V$ possesses a stratification with dual combinatorics, whose strata are canonically isomorphic to $\Omega_{V'}$ for all non-zero subspaces $V'$ of~$V$. We show that $Q_V$ is regular along all strata of codimension~$1$ and singular along all strata of codimension $\ge2$. In particular $Q_V$ as a whole is regular if and only if $r\le2$.

\medskip
The stratification has a natural description in terms of a modular interpretation of~$Q_V$, which is explained in Section~\ref{modonly}. Just as $P_V$ represents a certain functor of $\BF_q$-linear maps, the scheme $Q_V$ represents a functor of what we call \emph{reciprocal maps}. We do not know whether this somewhat strange concept has other uses. 

\medskip
The same goes for the natural morphisms $Q_V\to P_{V^*}$ and $P_{V^*}\to Q_V$ defined in Section~\ref{maps}, where $V^*$ denotes the vector space dual to~$V$. Their composites in both directions are a certain power of Frobenius; hence these morphisms are bijective and radicial. Whatever their deeper meaning, if any, they map strata to strata and thereby explain again why the combinatorics of the above mentioned stratifications correspond so literally.

\medskip
Finally recall that both $P_V$ and $Q_V$ are compactifications of the same affine variety~$\Omega_V$. In the last Section \ref{BV} we construct a third compactification $B_V$ which dominates both $P_V$ and $Q_V$ and which is smooth. We show that $B_V$ possesses a natural stratification indexed by flags of~$V$ and that the complement $B_V\setminus\Omega_V$ is a divisor with normal crossings. Thus $B_V$ constitutes a resolution of singularities of~$Q_V$ in the best possible sense.

\bigskip
The authors express their gratitude to Florian Breuer and Andrew Kresch for valuable comments.



\section{The ring $R_V$ and the variety $Q_V$}
\label{RV1}

Let $\BF_q$ be a finite field with $q$ elements, which we fix throughout the article. 
For any $\BF_q$-vector space $V$ we set
$$\circV \ :=\ V \setminus \{0\}.$$
For any non-zero finite dimensional $\BF_q$-vector space $V$ we define:
\begin{eqnarray}
\nonumber  S_V   & := & \text{the symmetric algebra of $V$ over $\BF_q$,} \\
\nonumber  K_V   & := & \mbox{the field of quotients of $S_V$,}\\
\nonumber  R_V & := & \mbox{the $\BF_q$-subalgebra of $K_V$ generated by $\frac{1}{v}$ for all $v \in \circV$,} \\
\nonumber  RS_V & := & \mbox{the $\BF_q$-subalgebra of $K_V$ generated by $R_V$ and $S_V$.}
\end{eqnarray}
Thus $RS_V$ is the localization of $S_V$ obtained by inverting all $v \in \circV$, and also the localization of $R_V$ obtained by inverting all $\frac{1}{v}$ for $v \in \circV$. Moreover, $K_V$ is also the field of quotients of~$R_V$. 

\medskip
The $\BF_q$-algebras $S_V$, $R_V$, and $RS_V$ are naturally $\BZ$-graded such that all $v\in\circV$ are homogeneous of degree $1$ and their reciprocals $\frac{1}{v}$ homogeneous of degree~$-1$. We indicate the homogenous parts of degree $n$ by $S_{V,n}$, $R_{V,n}$, and $RS_{V,n}$. Note that $S_{V,-n} = R_{V,n}=0$ for $n>0$, while $RS_V$ lives in all degrees.

\medskip
The generators of $R_V$ satisfy the following fundamental identities:
\begin{eqnarray}
\addtocounter{Thm}{1}\label{Relation1}
\frac{1}{\alpha v} &=& \rlap{\mbox{$\displaystyle\alpha^{-1}\cdot\frac{1}{v}$}}
\hphantom{\frac{1}{v+v'}\cdot\biggl(\frac{1}{v} + \frac{1}{v'}\biggr)}
\qquad\mbox{for all $v\in\circV$ and $\alpha\in\BF_q^\times$, and}\\[4pt]
\addtocounter{Thm}{1}\label{Relation2}
\qquad \frac{1}{v}\cdot\frac{1}{v'} &=& 
\frac{1}{v+v'}\cdot\biggl(\frac{1}{v} + \frac{1}{v'}\biggr)
\qquad\mbox{for all $v$, $v' \in \circV$ such that $v + v' \in\circV$.}
\end{eqnarray}
A useful reformulation of the second identity is:
\begin{eqnarray}
\addtocounter{Thm}{1}\label{Relation3}
\qquad \frac{1}{v}\cdot\frac{1}{v'} &=& 
\frac{1}{v-v'}\cdot\biggl(\frac{1}{v'} - \frac{1}{v}\biggr)
\qquad\mbox{for all $v$, $v' \in \circV$ such that $v-v' \in\circV$.}
\end{eqnarray}
We will use these identities to present $R_V$ by generators and relations. Let $A_V$ 
denote the polynomial ring over $\BF_q$ in the indeterminates $Y_v$ for all $v\in\circV$.  Let $\Fa_V \subset A_V$ be the homogeneous ideal generated by all elements of the form
\addtocounter{Thm}{1}
\begin{equation}
\label{RV1:generators}
\begin{cases}
\ Y_{\alpha v} - \alpha^{-1} Y_v
& \mbox{for all $v \in \circV$ and $\alpha \in \BF_q^{\times}$, and} \\[5pt]
\ Y_v Y_{v'} - Y_{v+v'} \cdot (Y_{v}+Y_{v'})
& \mbox{for all $v$, $v' \in \circV$ such that $v + v' \in\circV$.}
\end{cases}
\end{equation}
The identities (\ref{Relation1}) and (\ref{Relation2}) imply that $\Fa_V$ is contained in the kernel of the surjective $\BF_q$-algebra homomorphism $A_V \onto R_V$ defined by $Y_v\mapsto\frac{1}{v}$. We thus obtain a surjection
\addtocounter{Thm}{1}
\begin{equation}
\label{RV1:surjection}
A_V/\Fa_V \onto R_V.
\end{equation}
The following theorems are proved in Section \ref{RV2}:

\begin{Thm}
\label{RV1:presentation}
The homomorphism (\ref{RV1:surjection}) is an isomorphism.
\end{Thm}

\begin{Thm}
\label{RV1:CMnormal}
The ring $R_V$ is a Cohen-Macaulay normal integral domain.
\end{Thm}

\begin{Rem}
\label{RV1:UFD}
\rm If $\dim V=1$, the ring $R_V$ is isomorphic to a polynomial ring in one variable over~$\BF_q$, namely in $\frac{1}{v}$ for any $v\in\circV$. But if $\dim V\ge2$, it is not even factorial, because the identity (\ref{Relation2}) for linearly independent $v$, $v'$ gives two inequivalent factorizations whose factors are homogeneous of degree $-1$ and therefore indecomposable.
\end{Rem}

Next, let $|I|$ denote the cardinality of a set~$I$. For any integer $r\ge1$ we write $\{2,\ldots,r\}$ for the set of integers $i$ satisfying $2\le i\le r$, which is the empty set if $r=1$. Consider the polynomial
\addtocounter{Thm}{1}
\begin{equation}
\label{RV1:Hr}
h_r(T) := \sum_{I\subset\{2,\ldots,r\}} q^{\sum_{i\in I}(i-1)} \cdot \binom{T}{|I|} \ \in\ \BQ[T],
\end{equation}
where $I$ runs through all subsets of $\{2,\ldots,r\}$, including the empty set. Some initial cases are
\begin{eqnarray*}
h_1(T) &=& 1,\\
h_2(T) &=& 1+qT, \\
h_3(T) &=& 1+qT+q^2T+\smash{\frac{q^3}{2}}(T^2-T).
\end{eqnarray*}

\begin{Thm}
\label{RV1:hilbertfunction}
For $r:=\dim V$ and all $n\in\BZ$ the homogeneous part of $R_V$ of degree $-n$ has dimension
$$\dim R_{V,-n} = 
\begin{cases} h_r(n) & \text{if $n \ge 0$,} \\ 0 & \text{if $n<0$.} \end{cases}$$
\end{Thm}

Now we interpret the above rings as coordinate rings of algebraic varieties over~$\BF_q$. By construction $S_V$ is isomorphic to the polynomial ring $\BF_q[X_1,\ldots,X_r]$ with $r:=\dim V$. Thus $P_V := \Proj S_V$ is isomorphic to the standard projective space $\BP^{r-1}_{\BF_q}$. The localization $RS_V$ of $S_V$ defines an affine scheme $\Omega_V := \Spec RS_{V,0}$ that can be viewed as an open dense subscheme of~$P_V$. Under the identification with standard projective space it corresponds to $\BP^{r-1}_{\BF_q} \setminus ($union of all $\BF_q$-rational hyperplanes$)$.

\medskip
Usually the coordinate rings of projective algebraic varieties are assumed to be graded in degrees $\ge0$. But the construction of $\Proj R$ for a graded ring $R$ works equally for rings graded in degrees $\le0$. In this sense (or, if one prefers, with the grading inverted) we obtain another projective algebraic variety $Q_V := \Proj R_V$ over~$\BF_q$. The inclusion $R_V\into RS_V$ then also identifies $\Omega_V$ with an open dense subscheme of~$Q_V$. 
Thus we can view $Q_V$ as another interesting compactification of~$\Omega_V$ besides~$P_V$. Note that the presentation \ref{RV1:presentation} of $R_V$ describes $Q_V$ as the subvariety of $\BP^{q^r-2}_{\BF_q}$ determined by the ideal~$\Fa_V$. Theorem \ref{RV1:CMnormal} implies:

\begin{Thm}
\label{QV1:CMnormal}
The variety $Q_V$ is integral, Cohen-Macaulay, and projectively normal.
\end{Thm}



\section{Induction proofs}
\label{RV2}

The proofs of the theorems from Section~\ref{RV1} will use a basis of $V$ and repeated induction on $\dim V$. To facilitate this in a comprehensive setup we choose an infinite sequence of $\BF_q$-vector spaces $V_0 \subset V_1 \subset \ldots$ with $\dim V_r = r$ and prove the theorems for these. 

\medskip
We fix independent variables $X_1,X_2,\ldots$. For any integer $r\ge0$ we let $V_r$ denote the $\BF_q$-vector space with basis $X_1,\ldots,X_r$. For any $r\ge1$ the ring $R_{V_r}$ is contained in the rational function field $K_{V_r} = \BF_q(X_1,\ldots,X_r)$, and we have natural inclusions $0 = V_0 \subset V_1 \subset \ldots$ and $R_{V_1} \subset R_{V_2} \subset \ldots$ and $K_{V_1} \subset K_{V_2} \subset \nobreak\ldots$. For any $r \ge 1$ we define
\begin{eqnarray}
\addtocounter{Thm}{1}\label{RV2:FrDef}
f_r & := & \sum_{u \in V_{r-1}} \frac{1}{X_r + u}  \ \in \ R_{V_r}, \\
\addtocounter{Thm}{1}\label{RV2:DeltarDef}
\Delta_r & := & \bigl\{ 1 \bigr\} \cup 
           \biggl\{ \frac{1}{X_r + u} \,\biggm|\, u \in \circV_{r-1} \biggr\}
\ \subset \ R_{V_r},
\end{eqnarray}
and similarly
\begin{eqnarray}
\nonumber  \tilde{f_r} & := & \sum_{u \in V_{r-1}} Y_{X_r + u}  \ \in \ A_{V_r}, \\
\nonumber  \tilde{\Delta}_r & := & \bigl\{1\bigr\} \cup \bigl\{ Y_{X_r + u} \,\bigm|\, u \in \circV_{r-1} \bigr\} \ \subset \ A_{V_r}.
\end{eqnarray}
Then the surjection $A_{V_r} \onto R_{V_r}$ sends $\tilde{f_r}$ to $f_r$ and maps $\tilde{\Delta}_r$ bijectively to~$\Delta_r$. Note that we have natural inclusions $A_{V_1} \subset A_{V_2} \subset \ldots$ and $\Fa_{V_1} \subset \Fa_{V_2} \subset \ldots$.


\begin{Lem}
\label{RV2:decomposition}
For any $r\ge2$ we have
$$A_{V_r}\ =\ \Fa_{V_r} + \sum_{\tilde\delta \in \tilde{\Delta}_r} 
              A_{V_{r-1}}\bigl[\tilde{f_r}\bigr] \cdot \tilde\delta.$$
\end{Lem}

\begin{Proof}
We denote the right hand side of the equation by~$M$. Clearly $M$ is an $A_{V_{r-1}}\bigl[\tilde{f_r}\bigr]$-submodule of~$A_{V_r}$. 
We proceed in several steps.

\medskip
{\it (a) $1\in M$.} 

\medskip
This follows directly from $1 \in \tilde{\Delta}_r$.

\medskip
{\it (b) $Y_v \in M$ for all $v \in \circV_{r}$.}

\medskip
If $v\in\circV_{r-1}$ this follows from $Y_v \in A_{V_{r-1}}$ and $1 \in \tilde{\Delta}_r$.
Otherwise $v = \alpha \cdot (X_r + u)$ for some $u \in V_{r-1}$ and $\alpha  \in \BF_q^\times$. Then
$$Y_v = Y_{\alpha \cdot (X_r + u)}  \equiv \alpha^{-1} \cdot Y_{X_r + u} \bmod{\mathfrak{a}_{V_r}},$$ 
so it suffices to show that $Y_{X_r + u} \in M$. If $u \neq 0$, this follows from $Y_{X_r + u} \in \tilde{\Delta}_r$. If $u = 0$, then
$$Y_{X_r} = \tilde{f_r} - \sum_{w \in \circV_{r-1}} Y_{X_r +w},$$
where all summands of the right hand side lie in~$\tilde{\Delta}_r$; hence $Y_{X_r} \in M$, as desired.

\medskip
{\it (c) $Y_v Y_{v'} \in M$ for all $v, v' \in \circV_r$.}

\medskip
If $v\in\circV_{r-1}$, we have $Y_v Y_{v'} \in A_{V_{r-1}} M \subset M$ by~(b), and likewise if $v'\in\circV_{r-1}$. 
Otherwise $v = \alpha \cdot (X_r +u)$ and $v' = \alpha' \cdot (X_r + u')$ for some $u, u' \in V_{r-1}$ and $\alpha, \alpha' \in \BF_q^\times$. As in (b) we can reduce ourselves to the case that $\alpha=\alpha'=1$.
If $u \neq u'$, then
$$Y_{X_r +u} \cdot Y_{X_r + u'} \equiv Y_{u-u'} \cdot (Y_{X_r + u'} - Y_{X_r + u}) \bmod{\Fa_{V_r}}.$$
Here the right hand side lies in $A_{V_{r-1}}(M-M) \subset M$ by (b), as desired. If $u = u'$, then 
\begin{eqnarray}
\nonumber Y_{X_r +u} \cdot Y_{X_r + u} 
& = & \biggl(\tilde{f_r} - \sum_{\myatop{w \in V_{r-1}}{w \neq u}} Y_{X_r + w}\biggr) 
\cdot Y_{X_r + u}  \\
\nonumber & = & \tilde{f_r} \cdot Y_{X_r + u} - \sum_{\myatop{w \in V_{r-1}}{w \neq u}} Y_{X_r + w} \cdot Y_{X_r + u}.
\end{eqnarray}
Here the first summand lies in $A_{V_{r-1}}M \subset M$ by (b) and the remaining ones lie
in $M$ by the preceding case. Thus (c) follows in all cases.

\medskip
{\it (d) $M$ is an ideal of $A_{V_r}$.}

\medskip
Since $A_{V_r}$ is the $\BF_q$-algebra generated by $Y_v$ for all $v\in\circV_r$, it suffices to show that $Y_v M\subset M$. As $\Fa_{V_r}$ is already an ideal, it is enough to prove that $Y_v\tilde\Delta_r \subset M$. But this follows from (b) and~(c).

\medskip
Finally, assertions (a) and (d) together imply that $M =A_{V_r}$.
\end{Proof}


\begin{Lem}
\label{RV2:generation}
For any $r\ge2$ the ring $R_{V_r}$ is generated by $\Delta_r$ as a module over the subring $R_{V_{r-1}}[f_r]$.
\end{Lem}

\begin{Proof}
Direct consequence of the surjection (\ref{RV1:surjection}) and Lemma \ref{RV2:decomposition}.
\end{Proof}

\begin{Lem}
\label{RV2:transcendental}
For any $r\ge2$ the element $f_r \in K_{V_r}$ is transcendental over~$K_{V_{r-1}}$.
\end{Lem}

\begin{Proof}
Lemma \ref{RV2:generation} implies that the field extension $K_{V_{r-1}}(f_r) \subset K_{V_r}$ is finite. Since the transcendence degrees of $K_{V_{r-1}}$ and $K_{V_r}$ differ by~$1$, we conclude that $f_r$ must be transcendental over $K_{V_{r-1}}$.
\end{Proof}

\begin{Prop}
\label{RV2:free1}
For any $r\ge2$ the ring $R_{V_r}$ is a free module of rank $q^{r-1}$ with basis $\Delta_r$ over the subring $R_{V_{r-1}}[f_r]$.
\end{Prop}

\begin{Proof}
Let $W_r$ denote the group of automorphisms of $V_r$ that fix every element of $V_{r-1}$ and send $X_r$ to $X_r+u$ for some $u\in V_{r-1}$. This is a finite group of order $|V_{r-1}|=q^{r-1}$ that acts faithfully on~$V_r$. By functoriality it thus acts faithfully on $K_{V_r}$ and restricts to the identity on $K_{V_{r-1}}$. By (\ref{RV2:FrDef}) it also fixes~$f_r$; hence we obtain field inclusions
$$K_{V_{r-1}}(f_r) \subset K_{V_r}^{W_r} \subset K_{V_r}.$$
By Galois theory the extension $K_{V_r}^{W_r} \subset K_{V_r}$ has degree~$q^{r-1}$. But by Lemma \ref{RV2:generation}, the extension $K_{V_{r-1}}(f_r) \subset K_{V_r}$ has degree at most $|\Delta_r| =q^{r-1}$. Thus $K_{V_{r-1}}(f_r) = K_{V_r}^{W_r}$ and $\Delta_r$ is linearly independent over $K_{V_{r-1}}(f_r)$. In particular, $\Delta_r$ is linearly independent over $R_{V_{r-1}}[f_r]$, so the statement follows in conjunction with Lemma \ref{RV2:generation}.
\end{Proof}


\medskip
\begin{Proofof}{Theorem \ref{RV1:presentation}}
We must show that for all $r \ge 1$, the surjection $A_{V_r}/\Fa_{V_r}\onto R_{V_r}$
from (\ref{RV1:surjection}) is an isomorphism.
For $r=1$ this is a direct consequence of the identity (\ref{Relation1}).
Assume that $r\ge2$ and that the assertion holds for $r-1$. 
Take a new variable $T$ and consider the ring homomorphism $(A_{V_{r-1}}/ \Fa_{V_{r-1}})[T] \to A_{V_r}/\Fa_{V_r}$ induced by the inclusion $A_{V_{r-1}}\subset A_{V_r}$ and by $T\mapsto\tilde{f_r}$. By Lemma \ref{RV2:decomposition} this turns $A_{V_r}/\Fa_{V_r}$ into an $(A_{V_{r-1}}/ \Fa_{V_{r-1}})[T]$-module that is generated by the image of~$\tilde\Delta_r$.
On the other hand, by the induction hypothesis and Lemma \ref{RV2:transcendental} we have $(A_{V_{r-1}}/ \Fa_{V_{r-1}})[T] \cong R_{V_{r-1}}[T] \cong R_{V_{r-1}}[f_r]$. Thus Proposition \ref{RV2:free1} shows that $R_{V_r}$ becomes a free module with basis $\Delta_r$ over $(A_{V_{r-1}}/ \Fa_{V_{r-1}})[T]$.
Together we find that $A_{V_r}/\Fa_{V_r}\onto R_{V_r}$ is a surjective homomorphism of $(A_{V_{r-1}}/ \Fa_{V_{r-1}})[T]$-modules that sends a finite set of generators to a basis of the free module $R_{V_r}$. The homomorphism is therefore an isomorphism, as desired.
\end{Proofof}


\begin{Thm}
\label{RV2:free2}
\begin{itemize}
\item[(a)] The elements $f_1, \ldots, f_r \in R_{V_r}$ are algebraically independent over $\BF_q$. 
\item[(b)] The ring $R_{V_r}$ is a free module of rank $q^{\frac{r(r-1)}{2}}$ with basis $\Delta_1 \cdots \Delta_r$ over $\BF_q[f_1, \ldots, f_r]$.
\end{itemize}
\end{Thm}

\begin{Proof}
The case $r=1$ with $f_1=X_1$ is clear. The general case follows by induction on $r$ using Lemma \ref{RV2:transcendental} and Proposition \ref{RV2:free1}.
\end{Proof}

\medskip
\begin{Proofof}{Theorem \ref{RV1:CMnormal}}
By Theorem \ref{RV2:free2} the ring $R_{V_r}$ is free of finite rank over the polynomial ring $\BF_q[f_1, \ldots, f_r]$. Thus the elements $f_1, \ldots, f_r$ form a regular sequence in~$R_{V_r}$ of length equal to the Krull dimension of~$R_{V_r}$. The same then follows for the localization $R_{V_r,\Fm}$ of $R_{V_r}$ at the irrelevant maximal ideal $\Fm := \bigoplus_{n>0} R_{V_r,-n}$; hence this localization is Cohen-Macaulay. Using \cite[Cor.$\;$2.2.15]{Bruns-Herzog}
it follows that the graded ring $R_{V_r}$ itself is Cohen-Macaulay.

In particular $R_{V_r}$ satisfies Serre's condition (S${}_2$) (see \cite[p.$\;$63]{Bruns-Herzog}). By construction it is an integral domain. On the other hand, Theorem \ref{mod:singularlocus} below---whose proof does not depend on Theorem \ref{RV1:CMnormal}---implies that $\Proj R_{V_r}$ and therefore $R_{V_r}$ is regular in codimension one. Thus $R_{V_r}$ satisfies Serre's conditions (R${}_1$) and (S${}_2$) and is therefore normal (see \cite[Th.$\;$2.2.22]{Bruns-Herzog}), as desired.
\end{Proofof} 


\medskip
Next, for any $r \ge 1$ consider the subset of cardinality $q^{r-1}$
\addtocounter{Thm}{1}
\begin{equation}
\label{RV2:DeltaPrDef}
E_r \ :=\ \biggl\{ \frac{1}{X_r + u} \,\biggm|\, u \in V_{r-1} \biggr\}
\ \subset \ R_{V_r}.
\end{equation}
Note that it agrees with $\Delta_r$ except that $1$ has been replaced by $\frac{1}{X_r}$.
For any subset $I\subset\{1,\ldots,r\}$ we abbreviate 
\addtocounter{Thm}{1}
\begin{equation}
\label{RV2:EIDef}
E_I\ :=\ \prod_{i\in I}E_i.
\end{equation}

\begin{Lem}
\label{RV2:mydecomp}
For any $r\ge1$ we have
$$R_{V_r}\ =\ 
\bigoplus_{I\subset\{2,\ldots,r\}} \biggl(\;
\bigoplus_{e\in E_I}
\BF_q\bigl[f_1,f_i|_{i\in I}\bigr] \cdot e \biggr) .$$
\end{Lem}

\begin{Proof}
For $r=1$ the assertion follows from $R_{V_1} = \BF_q\bigl[\frac{1}{X_1}\bigr] = \BF_q[f_1]$. For $r\ge2$ Proposition \ref{RV2:free1} implies that 
$$R_{V_r} \ =\ \Bigl( R_{V_{r-1}} \oplus R_{V_{r-1}}[f_r] \cdot f_r \Bigr)
\oplus \bigoplus_{u\in \circV_{r-1}} R_{V_{r-1}}[f_r] \cdot\frac{1}{X_r+u}.$$
Using the definition (\ref{RV2:FrDef}) of $f_r$ we can rewrite this in the form
$$R_{V_r} \ =\ R_{V_{r-1}} \oplus \bigoplus_{u\in V_{r-1}} R_{V_{r-1}}[f_r] \cdot\frac{1}{X_r+u}.$$
From this the proposition follows by induction on~$r$.
\end{Proof}

\medskip
Let $U_r$ denote the group of automorphisms of $V_r$ that send each $X_i$ to $X_i+u_i$ for some $u_i\in V_{i-1}$. It corresponds to the group of upper-triangular matrices in $\GL_r(\BF_q)$ with all diagonal entries equal to~$1$. Let $U_r$ act on $R_{V_r}$ by functoriality.

\begin{Thm}
\label{RV2:invariants}
The ring of invariants is
$R_{V_r}^{U_r} = \BF_q[f_1, \ldots, f_r].$
\end{Thm}

\begin{Proof}
By construction $U_r$ fixes each $f_i$; this implies `$\supset$'. By construction $U_r$ also
acts transitively on $E_I$ for each subset $I\subset\{2,\ldots,r\}$. Therefore it stabilizes the corresponding inner sum in Lemma \ref{RV2:mydecomp} and permutes its basis transitively. The $U_r$-invariants in that part thus form the $\BF_q[f_1,f_i|_{i\in I}]$-submodule generated by
$\sum_{e\in E_I}e = \prod_{i\in I}f_i$. Together this proves
 `$\subset$', and we are done.
\end{Proof}



\begin{Prop}
\label{RV2:HomRep}
For any $r\ge1$ and $n\ge0$, there is an isomorphism of representations of $U_r$ over~$\BF_q$:
$$R_{V_r,-n}\ \cong\ 
\bigoplus_{I\subset\{2,\ldots,r\}} \biggl(\;
\bigoplus_{e\in E_I}
\BF_q\cdot e \biggr)^{\textstyle\binom{n}{|I|}} .$$
\end{Prop}

\begin{Proof}
Since each $f_i$ and each element of $E_i$ is homogeneous of degree $-1$, the decomposition in Lemma \ref{RV2:mydecomp} is graded and each $e\in E_I$ is homogeneous of degree~$-|I|$. Thus it suffices to show that for any $I\subset\{2,\ldots,r\}$, the homogeneous part of degree $-n+|I|$ of $\BF_q\bigl[f_1,f_i|_{i\in I}\bigr]$ has dimension $\binom{n}{|I|}$. Set $d:=n-|I|$ and $k:=|I|+1$, so that $k\ge1$ and $d = n+1-k\ge1-k$. Then $\BF_q\bigl[f_1,f_i|_{i\in I}\bigr]$ is isomorphic to a polynomial ring in $k$ variables that are homogeneous of degree~$-1$. Thus after inverting degrees we need to know that the homogeneous part of degree $d$ in a polynomial ring in $k\ge1$ variables has dimension $\binom{d+k-1}{k-1}$ whenever $d\ge1-k$. This is the well-known formula if $d\ge0$, and for $0>d\ge1-k$ it holds because both sides are zero.
\end{Proof}

\medskip
\begin{Proofof}{Theorem \ref{RV1:hilbertfunction}}
Forgetting the action of $U_r$ in Proposition \ref{RV2:HomRep}, we find that 
$$\dim R_{V_r,-n} \ = 
\sum_{I\subset\{2,\ldots,r\}} 
\bigl|E_I\bigr| \cdot {\textstyle\binom{n}{|I|}} \ =
\sum_{I\subset\{2,\ldots,r\}} \Bigl(\prod_{i\in I}q^{i-1}\Bigr) 
\cdot{\textstyle\binom{n}{|I|}}
\ =\ h_r(n)$$
whenever $n\ge0$. Since $R_{V_r,n}=0$ for $n>0$, the theorem follows.
\end{Proofof}



\section{Rings of invariants and quotient varieties}
\label{invar}

Let $V$ be an $\BF_q$-vector space of finite dimension $r\ge1$. Let $G:=\Aut_{\BF_q}(V)$ be its automorphism group and $G'\subset G$ the subgroup of automorphisms of determinant~$1$. Then of course $G\cong\GL_r(\BF_q)$ and $G'\cong \SL_r(\BF_q)$. Let $U$ be a maximal unipotent subgroup of~$G$; in a suitable basis of $V$ it corresponds to the group of upper-triangular matrices with all diagonal entries equal to~$1$. By functoriality these groups act on $R_V$ and~$S_V$.

\begin{Thm}
\label{invar:Rings}
The respective subrings of invariants are generated by algebraically independent homogeneous elements as follows:
\begin{itemize}
\item[(a)] 
$R_V^U = \BF_q[f_1,\ldots,f_r]$ 
with all degrees $-1$.
\item[(b)] 
$S_V^U = \BF_q[g_1,\ldots,g_r]$ 
with respective degrees $1,\,q,\ldots,\,q^{r-1}$.
\item[(c)] 
$R_V^G = \BF_q[h_1,\ldots,h_r]$ 
with respective degrees $1{-}q,\,1{-}q^2,\ldots,\,1{-}q^r$.
\item[(d)] 
$S_V^G = \BF_q[k_0,\ldots,k_{r-1}]$ 
with respective degrees $q^r{-}1,\,q^r{-}q,\ldots,\,q^r{-}q^{r-1}$.
\item[(e)] 
$R_V^{G'} = \BF_q[h_1,\ldots,h_{r-1},h'_r]$ 
with respective degrees $1{-}q,\ldots,\,1{-}q^{r-1},\,\smash{\frac{1-q^r}{q-1}}$.
\item[(f)] 
$S_V^{G'} = \BF_q[k'_0,k_1,\ldots,k_{r-1}]$ 
with respective degrees $\frac{q^r-1}{q-1},\,q^r{-}q,\ldots,\,q^r{-}q^{r-1}$.
\end{itemize}
\end{Thm}

\begin{Proof}
We use the notations from Section~\ref{RV2} and identify $V=V_r$ and $U=U_r$. Then (a) is just the combination of Theorems \ref{RV2:free2}~(a) and~\ref{RV2:invariants}. Next $S_{V_r}$ is the polynomial ring $\BF_q[X_1,\ldots,X_r]$. For any $1\le i\le r$ its element $g_i := \prod_{u\in V_{i-1}}(X_i+u)$ is invariant under~$U_r$ and homogeneous of degree~$q^{i-1}$. By induction on $r$ one easily shows that $S_{V_r}$ is a free module over the subring $\BF_q[g_1,\ldots,g_r]$ with basis 
$\{ X_1^{\nu_1}\cdots X_r^{\nu_r} \mid \forall i:\; 0\le\nu_i<q^i-1 \}.$
In particular it is free of rank $q^{\frac{r(r-1)}{2}} = |U_r|$. Since by Galois theory the quotient field extension $K_{V_r}/K_{V_r}^{U_r}$ also has degree $|U_r|$, we deduce that $K_{V_r}^{U_r}$ is the quotient field of $\BF_q[g_1,\ldots,g_r]$. But this ring is isomorphic to a polynomial ring and hence integrally closed in its quotient field. Thus (b) follows.

\medskip
Next take an auxiliary variable $T$. A classical theorem of Dickson \cite[Th.$\;$1.2]{Wilkerson} states that 
$$k(T)\ :=\ \prod_{v\in V} (T-v)\ =\ T^{q^r} + \sum_{i=0}^{r-1} k_i\, T^{q^i}$$
for algebraically independent elements $k_i \in S_V$ and that $S_V^G = \BF_q[k_0,\ldots,k_{r-1}]$.
Since the defining equation is jointly homogeneous in $T$ and $v\in V$, we find that $k_i$ is homogeneous of degree $q^r-q^i$. This proves~(d). For the proof of (c), which follows a suggestion of Florian Breuer, we calculate
$$h(T)\ :=\ \prod_{v\in\circV} {\bigl(\textstyle T-\frac{1}{v}\bigr)}
\ =\ \prod_{v\in\circV} \frac{T}{-v}\cdot\Bigl(\frac{1}{T}-v\Bigr)
\ =\ \frac{T^{q^r}}{k_0} \cdot k\bigl(\textstyle\frac{1}{T}\bigr).$$
It follows that
$$h(T)\ =\ T^{q^r-1} + \sum_{i=1}^r h_i T^{q^r-q^i}$$ 
with coefficients $h_i:= \frac{k_i}{k_0}$ for $1\le i\le r-1$ and $h_r := \frac{1}{k_0}$, which are homogeneous of the indicated degrees. By construction $h(T)$ has coefficients in~$R_V^G$; hence $\BF_q[h_1,\ldots,h_r] \subset R_V^G$. Moreover this ring extension is integral, because all generators $\frac{1}{v}$ of $R_V$ are zeroes of $h(T)$. On the other hand, the form of the $h_i$ implies that
$$K_V^G \ =\ \Quot(S_V^G) \ \stackrel{(d)}{=}\ \BF_q(k_0,\ldots,k_{r-1})
\ \stackrel{!}{=}\ \BF_q(h_1,\ldots,h_r) \ \subset \ \Quot(R_V^G) \ =\ K_V^G.$$
Since $\BF_q[h_1,\ldots,h_r]$ is integrally closed in its field of quotients, it is therefore equal to~$R_V^G$, proving~(c).

\medskip
For (f) let $k'_0 \in S_V$ denote the product of one non-zero element from every $1$-dimen\-sional subspace of~$V$ (it does not matter which). This is a homogeneous polynomial of degree $\frac{q^r-1}{q-1}$. Dickson \cite[Th.$\;$3.1]{Wilkerson} also proved that  $S_V^{G'} = \BF_q[k'_0,k_1,\ldots,k_{r-1}]$, whence~(f). Finally, this implies (e) in the same manner as (d) implies (c).
\end{Proof}


Next we recall the notion of a weighted projective space. For a general introduction to these see \cite{BeltramettiRobbiano}. Consider a polynomial ring $R=k[T_1,\ldots,T_r]$ over a field~$k$. Take positive integers $d_1,\ldots,d_r$ and endow $R$ with the unique grading for which each $T_i$ is homogeneous of degree~$d_i$. Then $\Proj R$ is called a \emph{weighted projective space of weights $d_1,\ldots,d_r$}. It is a normal projective algebraic variety.
Multiplying all $d_i$ by a fixed positive integer does not change $\Proj R$. The weighted projective space of weights $1,\ldots,1$ is just the usual projective space $\BP^{r-1}_k$.

\medskip
We are interested in the projective algebraic varieties with coordinate rings $S_V$ and $R_V$. Since by definition $R_V$ is graded in degrees $\le0$, we temporarily view it as graded in degrees $\ge0$ by multiplying all degrees by $-1$, and can then define
\begin{eqnarray*}
P_V &\!\!:=\!\!& \Proj S_V \ \cong\ \BP^{r-1}_k,\\
Q_V &\!\!:=\!\!& \Proj R_V.
\end{eqnarray*}
Everywhere else we will keep the previous grading on $R_V$. By functoriality the groups $U\subset G\supset G'$ act on $P_V$ and $Q_V$, and applying $\Proj$ to the respective subrings of invariants yields the associated quotient varieties. Thus Theorems \ref{invar:Rings} and \ref{RV2:free2} (b) imply:

\begin{Thm}
\label{invar:Varieties}
\begin{itemize}
\item[(a)] 
$Q_V/U \cong \BP^{r-1}_{\BF_q}$ and the projection $Q_V\onto Q_V/U$ is finite and flat of degree $q^{\frac{r(r-1)}{2}}$.
\end{itemize}
\medskip
The other quotients are weighted projective spaces of the following weights:
\begin{itemize}
\item[(b)] 
$P_V/U$ has weights $1,\,q,\ldots,\,q^{r-1}$.
\item[(c)] 
$Q_V/G$ has weights $q{-}1,\,q^2{-}1,\ldots,\,q^r{-}1$.
\item[(d)] 
$P_V/G$ has weights $q^r{-}1,\,q^r{-}q,\ldots,\,q^r{-}q^{r-1}$.
\item[(e)] 
$Q_V/G'$ has weights $q{-}1,\ldots,\,q^{r-1}{-}1,\,\smash{\frac{q^r-1}{q-1}}$.
\item[(f)]
$P_V/G'$ has weights $\frac{q^r-1}{q-1},\,q^r{-}q,\ldots,\,q^r{-}q^{r-1}$.
\end{itemize}
\end{Thm}

In particular $Q_V/U$ is always regular. The other five quotients are normal algebraic varieties of dimension $r-1$; hence they are regular if $r\le2$. By contrast:

\begin{Prop}
\label{invar:Regularity} 
The quotients in \ref{invar:Varieties} (b--f) are singular if $r\ge3$.
\end{Prop}

\begin{Proof}
A weighted projective space of weights $d_1,\ldots,d_r$ is regular if and only if for every prime number~$\ell$, the maximum of $\ord_\ell(d_1),\ldots,\ord_\ell(d_r)$ is attained at least $r-1$ times
(combine \cite[Lem.$\,$3.C.4, Prop.$\,$3.C.5, Prop.$\,$4.A.6$\,$(c)]{BeltramettiRobbiano}).
This criterion fails in the cases (b), (d), and~(f) for the prime $\ell\,|\,q$.
In the case (c) we can first divide all weights by the common factor $q-1$. Then the first weight is~$1$, and the criterion would require that all other weights are equal, which is clearly not the case. 
In the case (e) the criterion would also have to hold for the weights $q{-}1,\ldots,\,q^{r-1}{-}1$, where it fails by (c) if $r\ge4$. The remaining case $r=3$ of (e) is left to the careful reader.
\end{Proof}



\section{Invariants under unipotent subgroups}
\label{URep}

As before let $r:=\dim V \ge1$ and $G:=\Aut_{\BF_q}(V)\cong\GL_r(\BF_q)$. Let $H\subset G$ be a unipotent subgroup, or equivalently a $p$-subgroup, where $p$ is the characteristic of~$\BF_q$. In this section we study the subring of invariants~$R_V^H$.

\medskip
Let $U$ be a maximal unipotent subgroup of~$G$. We choose a basis of $V$ such that $U$ corresponds to the group of upper triangular matrices with all diagonal entries equal to~$1$. Then for any integer $1\le s\le r$ we let $P_s\subset G$ denote the subgroup corresponding to matrices of block triangular form
$$\begin{pmatrix} * & * \\ 0 & * \end{pmatrix}$$
where the upper left block has size $s\times s$ and the lower right has size $(r-s)\times(r-s)$. This is a maximal parabolic subgroup of $G$ if  $1\le s<r$, and equal to $G$ if $s=r$. Let $L_s$ denote the normal subgroup of $P_s$ corresponding to matrices of block triangular form
$$\begin{pmatrix} 1_s & * \\ 0 & * \end{pmatrix}$$
where $1_s$ denotes the identity matrix of size $s\times s$. Then $UL_s/L_s$ is a maximal unipotent subgroup of $P_s/L_s\cong\GL_s(\BF_q)$; hence the well-known formula $|\GL_s(\BF_q)| = \prod_{i=1}^s (q^s-q^{s-i})$ implies that $[P_s:UL_s] = 
\prod_{i=1}^s (q^i-1)$. Our aim is to prove the following result:

\begin{Thm}
\label{URep:Hilbert}
For any unipotent subgroup $H\subset G$ and any integer $n\ge0$ we have
$$\dim R^H_{V,-n} \ =\ \sum_{s=1}^r\ \frac{\,|H\backslash G/L_s|\,}{[P_s:UL_s]} \cdot \binom{n-1}{s-1}.$$
\end{Thm}

For the proof we use the notations from Section~\ref{RV2} and identify $V=V_r$ and $U=U_r$. Note that both sides of the desired equality remain unchanged on conjugating $H$ by an element of~$G$. Thus without loss of generality we may assume that $H\subset U$. One part of the calculation then involves describing the left action of $U$ on $G/L_s$. Let $E_i$ and $E_I$ be as in (\ref{RV2:DeltaPrDef}) and (\ref{RV2:EIDef}).

\begin{Lem}
\label{URep:Lem}
For any $1\le s\le r$ the set $G/L_s$ with the left action of $U$ is isomorphic to the disjoint union
of $[P_s:UL_s]$ copies of the sets $E_I$ for all $I\subset\{1,\ldots,r\}$ with $|I|=s$.
\end{Lem}

\begin{Proof}
Let $W$ denote the group of permutation matrices in~$G$, which we can identify with the symmetric group~$S_r$. Then $W$ is the Weyl group of $G$ and $W_s := W\cap P_s$ the Weyl group of~$P_s$, both with respect to the diagonal torus. By Bruhat we thus have the disjoint decomposition
$$G/L_s \ =\ \Bigl( \bigsqcup_{w\in W/W_s} UwP_s\Bigr) \!\Bigm/\!\! L_s\ =\ 
\bigsqcup_{w\in W/W_s} UwP_s/L_s.$$
Since $L_s$ is a normal subgroup of~$P_s$, the group $P_s$ still acts by right translation on $UwP_s/L_s$. This action commutes with left translation by~$U$, and the two actions together are transitive. Therefore all $U$-orbits in $UwP_s/L_s$ are isomorphic. 

\medskip
To determine their number observe that there is a bijection
$$(w^{-1}Uw\cap P_s)\backslash P_s/L_s\ \stackrel{\sim}{\longto}\ U\backslash UwP_s/L_s,\ \ 
[p] \mapsto [wp].$$
Since $L_s$ is a normal subgroup of~$P_s$, the left hand side is also the set of right cosets of the subgroup $(w^{-1}Uw\cap P_s)L_s$ of~$P_s$. We claim that $(w^{-1}Uw\cap P_s)L_s$ is conjugate to $UL_s$ under~$P_s$. Indeed, to prove this we can replace $w$ by any element of the coset $wW_s$. We can thus assume that $wi<wj$ for all $1\le i<j\le s$. Then 
$$(w^{-1}Uw) \cap \begin{pmatrix} * & 0 \\ 0 & 1_{r-s} \end{pmatrix}
\ =\ U \cap \begin{pmatrix} * & 0 \\ 0 & 1_{r-s} \end{pmatrix},$$
which implies that $(w^{-1}Uw\cap P_s)L_s = UL_s$, as desired. Using the claim we find that the number of $U$-orbits in $UwP_s/L_s$ is
$$\bigl|U\backslash UwP_s/L_s\bigr|\ =\ 
\bigl|(w^{-1}Uw\cap P_s)\backslash P_s/L_s\bigr|\ =\ 
\bigl[P_s: (w^{-1}Uw\cap P_s)L_s\bigr]\ =\ 
\bigl[P_s: UL_s\bigr].$$
The structure of a single $U$-orbit in $UwP_s/L_s$ is determined by the bijection
$$U / (U\cap wL_sw^{-1})\ \stackrel{\sim}{\longto}\ UwL_s/L_s\ \subset\ UwP_s/L_s,\ \ 
[u] \mapsto [uw].$$
Here $wL_sw^{-1}$ is the group of matrices in~$G$ whose $i$-th column coincides with that of the identity matrix for all $i\in I_w :=\{wi\mid 1\le i\le s\}$. Thus $U\cap wL_sw^{-1}$ is precisely the stabilizer in $U$ of the element 
$$ \prod_{i\in I_w} \frac{1}{X_i} \ \in\ \prod_{i\in I_w} E_i \ =\ E_{I_w}.$$
Since $U$ acts transitively on~$E_{I_w}$, we deduce that each $U$-orbit in $UwP_s/L_s$ is isomorphic to~$E_{I_w}$.

\medskip
Finally, the map $w\mapsto I_w$ induces a bijection from $W/W_s$ to the set of subsets $I\subset\{1,\ldots,r\}$ with $|I|=s$. By combining everything the lemma follows.
\end{Proof}

\begin{Proofof}{Theorem \ref{URep:Hilbert}}
Using Lemma \ref{URep:Lem} we can rewrite the right hand side of the desired equality as
$$\sum_{s=1}^r\ \frac{\,|H\backslash G/L_s|\,}{[P_s:UL_s]} \cdot \binom{n-1}{s-1}
\ = \sum_{\emptyset\not=I\subset\{1,\ldots,r\}} 
\bigl|H\backslash E_I\bigr|
 \cdot \binom{\,n-1}{|I|-1}.$$
Since $|E_1|=1$, the summand associated to $I=\{1\}$ is identically equal to~$1$. For the other terms the subset $I' := I\setminus\{1\}$ is again non-empty, and $|E_1|=1$ implies that 
$$\bigl|H\backslash E_I\bigr| \ =\ \bigl|H\backslash E_{I'}\bigr|.$$
Any non-empty subset $I'\subset\{2,\ldots,r\}$ arises in this way from the two subsets $I=I'$ and $I=\{1\}\sqcup I'$, and combining the corresponding terms the total sum becomes 
$$1+\!\! \sum_{\emptyset\not=I'\subset\{2,\ldots,r\}} 
\bigl|H\backslash E_{I'}\bigr|
 \cdot \biggl[ \binom{\,n-1}{|I'|-1} + \binom{n-1}{|I'|} \biggr].$$
By standard identities of binomial coefficients and the fact that $|E_\emptyset|=1$ this simplifies to
$$\sum_{I'\subset\{2,\ldots,r\}} 
\bigl|H\backslash E_{I'}\bigr|
\cdot \binom{n}{|I'|}.$$
At last, this is equal to $\dim R^H_{V,-n}$ by Proposition \ref{RV2:HomRep}.
\end{Proofof}



\section{Cohomology}
\label{cohom}

As before let $r=\dim V \ge1$. Let $\CO(1)$ denote the standard twisting sheaf on $P_V\cong\BP_{\BF_q}^{r-1}$ and recall that for all $i$, $n\in\BZ$ we have
\addtocounter{Thm}{1}
\begin{equation}
\label{cohom:PVCoh}
\dim H^i(P_V,\CO(n)) \ =\ 
\begin{cases} \binom{r-1+n}{r-1} & \text{if $i=0$ and $n\ge 0$,} \\ 
              (-1)^{r-1}\cdot \binom{r-1+n}{r-1} & \text{if $i=r-1$ and $n\le-r$,} \\ 
              0 & \text{otherwise.} \end{cases}
\end{equation}
Let $\CO_{Q_V}(1)$ denote the ample invertible sheaf on $Q_V$ corresponding to the graded $R_V$-module $R_V$ shifted in degrees by~$1$. Theorem \ref{invar:Rings} (a) directly implies:

\begin{Prop}
\label{cohom:O1}
$\CO_{Q_V}(1)$ is the pullback of the standard twisting sheaf $\CO(1)$ under the projection $\pi: Q_V\onto Q_V/U \cong \BP^{r-1}_{\BF_q}$.
\end{Prop}


\begin{Thm}
\label{cohom:QVCoh}
Let $h_r(T)$ be the polynomial from Theorem \ref{RV1:hilbertfunction} with $r=\dim V$. Then for all $i$, $n\in\BZ$ we have
$$\dim H^i(Q_V,\CO_{Q_V}(n)) \ =\ 
\begin{cases} h_r(n) & \text{if $i=0$ and $n\ge 0$,} \\ 
              (-1)^{r-1}\cdot h_r(n) & \text{if $i=r-1$ and $n<0$,} \\ 
              0 & \text{otherwise.} \end{cases}$$
\end{Thm}

\begin{Proof}
For $r=1$ both sides are always~$1$. So assume that $r\ge2$. Recall from (\ref{RV2:DeltarDef}) that $\Delta_1=\{1\}$ and that all other $\Delta_i$ consist of homogeneous elements of degree $0$ and~$-1$. Thus $\Delta_1 \cdots \Delta_r$ consists of some number $a_{r,s}$ of elements of degree $-s$ for all $0\le s\le r-1$. From Theorem \ref{RV2:free2} (b) we deduce that $\pi_*\CO_{Q_V}$ is isomorphic to the direct sum of $a_{r,s}$ copies of $\CO(-s)$ for all $0\le s\le r-1$. (The reader may check that this takes the sign convention for the grading on $R_V$ correctly into account.) In view of Proposition \ref{cohom:O1} we find that
\begin{eqnarray*}
    \dim H^i(Q_V,\CO_{Q_V}(n)) 
&=& \dim H^i(Q_V,\pi^*\CO(n)) \\
&=& \dim H^i(\BP^{r-1}_{\BF_q},(\pi_*\CO_{Q_V})(n)) \\
&=& \sum_{0\le s\le r-1} a_{r,s} \cdot \dim H^i(\BP^{r-1}_{\BF_q},\CO(n-s)).
\end{eqnarray*}
By (\ref{cohom:PVCoh}) this is zero unless $i=0$ or $r{-}1$, and for $i=0$ it is
$$\sum_{\myatop{0\le s\le r-1}{s\le n}} a_{r,s} \cdot \binom{r-1+n-s}{r-1}.$$
This again is zero unless $n\ge0$, and for such $n$ the sum can be extended over all $0\le s\le r{-}1$, because the binomial coefficient vanishes whenever $-(r{-}1)\le n{-}s<0$. Thus with 
$$k_r(T)\ :=\ \sum_{0\le s\le r-1} a_{r,s} \cdot \binom{r-1+T-s}{r-1}$$
we have $\dim H^0(Q_V,\CO_{Q_V}(n)) = k_r(n)$ for all $n\ge0$. 

\medskip
On the other hand, we have $\dim H^0(Q_V,\CO_{Q_V}(n)) = \dim R_{V,-n}$ for all $n\gg0$ by the general theory of projective coordinate rings \cite[Ch.$\,$II Exc.$\,$5.9]{Hartshorne}, which in turn is equal to $h_r(n)$ by Theorem~\ref{RV1:hilbertfunction}. Together we deduce that $k_r(n) = h_r(n)$ for all $n\gg0$. Since both sides are polynomials in~$n$, it follows that $k_r(T) = h_r(T)$.

\medskip
This now implies that $\dim H^0(Q_V,\CO_{Q_V}(n)) = h_r(n)$ for all $n\ge0$. A similar consideration as above shows that $\dim H^{r-1}(Q_V,\CO_{Q_V}(n))$ vanishes for all $n\ge0$. The formula for $n<0$ follows from the fact that the Euler characteristic is a polynomial in~$n$.
\end{Proof}

\begin{Cor}
\label{cohom:QVCohIsom}
For all $n$ the natural map $R_{V,-n} \longto H^0(Q_V,\CO_{Q_V}(n))$ is an isomorphism.
\end{Cor}



\section{Dualizing sheaf}
\label{dual}

As before let $r=\dim V \ge1$. Let $I_V \subset R_V$ be the homogeneous ideal generated by the elements $\frac{1}{v_0\cdots v_r}$ for all $v_0,\ldots,v_r\in \circV$, any $r$ of which are linearly independent. Let $\CI_V \subset \CO_{Q_V}$ be the ideal sheaf associated to~$I_V$, and let $\CI_V(r)$ be its $r$-fold twist by the ample invertible sheaf $\CO_{Q_V}(1)$ from Section~\ref{cohom}. The aim of this section is to prove:

\begin{Thm}
\label{dual:DualQV}
$\CI_V(r)$ is a dualizing sheaf on~$Q_V$.
\end{Thm}

We use the notations from Section \ref{RV2} and take $V=V_r$. We begin by describing a nice basis of~$I_{V_r}$. For any $r\ge1$ set
\begin{equation}
\addtocounter{Thm}{1}\label{dual:ErDef}
\hat\Delta_r\ :=\ \biggl\{ \frac{f_r}{X_r} \biggr\} \cup 
           \biggl\{ \frac{1}{X_r + u} - \frac{1}{X_r} \;\biggm|\; u \in \circV_{r-1} \biggr\}
\ \subset\ R_{V_r}.
\end{equation}
Let $J_r$ denote the $R_{V_{r-1}}[f_r]$-submodule of $R_{V_r}$ that is generated by~$\hat\Delta_r$. 

\begin{Lem}
\label{dual:JrIdeal}
$J_r$ is an ideal of $R_{V_r}$.
\end{Lem}

\begin{Proof}
It suffices to show that $\frac{1}{v}\cdot J_r\subset J_r$ for all generators $\frac{1}{v}$ of $R_{V_r}$ with $v\in\circV_r$. That is obvious for $v\in\circV_{r-1}$. By the identity (\ref{Relation1}) the other cases reduce to $\frac{1}{X_r+w}\cdot J_r\subset J_r$ for all $w\in V_{r-1}$. On the generators of $J_r$ this amounts to the relations
\begin{itemize}
\item[(a)] $\displaystyle\frac{1}{X_r+w}\cdot\biggl(\frac{1}{X_r+u}-\frac{1}{X_r}\biggr) \in J_r$
for all $w\in V_{r-1}$ and $u\in\circV_{r-1}$, and
\item[(b)] $\displaystyle\frac{1}{X_r+w}\cdot\frac{f_r}{X_r} \in J_r$ for all $w\in V_{r-1}$.
\end{itemize}
If $w\not=0$, $u$, by the identities (\ref{Relation3}) the element in (a) is
\begin{eqnarray*} 
&=& \frac{1}{w-u}\cdot\biggl(\frac{1}{X_r+u}-\frac{1}{X_r}\biggr)
+ \biggl(\frac{1}{w}-\frac{1}{w-u}\biggr)\cdot\biggl(\frac{1}{X_r+w}-\frac{1}{X_r}\biggr).
\end{eqnarray*} 
This is an $R_{V_{r-1}}$-linear combination of elements of~$\hat\Delta_r$; hence it lies in~$J_r$. If $w=0$, using the definition \ref{RV2:FrDef} of $f_r$ the element comes out as
\begin{eqnarray*} 
&=& \frac{1}{X_r} \cdot \frac{1}{X_r+u}
  - \frac{1}{X_r} \cdot \biggl( f_r - \sum_{v\in\circV_{r-1}} \frac{1}{X_r+v} \biggr) \\ 
&\!\!\stackrel{(\ref{Relation3})}{=}\!\!
& \frac{-1}{u}\cdot\biggl(\frac{1}{X_r+u}-\frac{1}{X_r}\biggr)
  - \frac{f_r}{X_r}
  - \sum_{v\in\circV_{r-1}} \frac{1}{v} \cdot \biggl(\frac{1}{X_r+v}-\frac{1}{X_r}\biggr),
\end{eqnarray*}
where the right hand side again lies in~$J_r$. If $w=u$, the element is
\begin{eqnarray*} 
&\!\!\stackrel{(\ref{RV2:FrDef})}{=}\!\!
& \biggl( f_r - \sum_{\myatop{v \in V_{r-1}}{w \neq u}} \frac{1}{X_r+v} \biggr)
  \cdot \biggl(\frac{1}{X_r+u}-\frac{1}{X_r}\biggr) \\
&=& f_r \cdot \biggl(\frac{1}{X_r+u}-\frac{1}{X_r}\biggr)
- \sum_{\myatop{v \in V_{r-1}}{w \neq u}}
\frac{1}{X_r+v} \cdot \biggl(\frac{1}{X_r+u}-\frac{1}{X_r}\biggr).
\end{eqnarray*}
Here the first summand lies in $f_r\cdot \hat\Delta_r$ and the remaining ones lie in $J_r$ by the preceding cases. This finishes the proof of~(a). For $w\not=0$ the element in (b) is
$$\stackrel{(\ref{Relation3})}{=}
\ \frac{-f_r}{w}\cdot\biggl(\frac{1}{X_r+w}-\frac{1}{X_r}\biggr)$$
and hence in~$J_r$. For $w=0$ it is 
\begin{eqnarray*}
&\!\!\stackrel{(\ref{RV2:FrDef})}{=}\!\!
& \biggl( f_r - \sum_{v\in\circV_{r-1}} \frac{1}{X_r+v} \biggr) \cdot \frac{f_r}{X_r} \\
&=& f_r\cdot\frac{f_r}{X_r}
- \sum_{v\in\circV_{r-1}} \frac{1}{X_r+v} \cdot \frac{f_r}{X_r}.
\end{eqnarray*}
Here the first summand lies in $f_r\cdot \hat\Delta_r$ and the remaining ones lie in $J_r$ by the preceding case. This finishes the proof of~(b).
\end{Proof}


\begin{Lem}
\label{dual:IrSubset}
$I_{V_{r-1}}\cdot J_r \subset I_{V_r}$ for any $r\ge2$. 
\end{Lem}

\begin{Proof}
Looking at generators, this amounts to showing that for all $v_0,\ldots,v_{r-1}\in \circV_{r-1}$, any $r-1$ of which are linearly independent, we have
\begin{itemize}
\item[(a)] $\displaystyle\frac{1}{v_0\cdots v_{r-1}} \cdot \biggl(\frac{1}{X_r+u}-\frac{1}{X_r}\biggr) \in I_{V_r}$
for all $u\in\circV_{r-1}$, and
\item[(b)] $\displaystyle\frac{1}{v_0\cdots v_{r-1}} \cdot \frac{f_r}{X_r} \in I_{V_r}$.
\end{itemize}
To prove (a) we keep $v_0,\ldots,v_{r-1}$ fixed and vary~$u$. Note that (a) trivially holds for $u=0$. For all $u\in V_{r-1}$ and $\alpha\in\BF_q^\times$ the identity (\ref{Relation3}) implies that 
\begin{eqnarray*}
\qquad\qquad\qquad\qquad\qquad\qquad\qquad\qquad\qquad && 
\llap{$\displaystyle\frac{1}{v_0\cdots v_{r-1}} \cdot \biggl(\frac{1}{X_r+u+\alpha v_i}-\frac{1}{X_r}\biggr)
- \frac{1}{v_0\cdots v_{r-1}}$} \cdot \biggl(\frac{1}{X_r+u}-\frac{1}{X_r}\biggr) \\
&=& \frac{1}{v_0\cdots v_{r-1}} \cdot \biggl(\frac{1}{X_r+u+\alpha v_0}-\frac{1}{X_r+u}\biggr) \\
&=& \frac{-\alpha}{v_1\cdots v_{r-1} \cdot (X_r+u+\alpha v_0)\cdot(X_r+u)}.
\end{eqnarray*}
Here any $r$ of the factors in the denominator on the right hand side are linearly independent; hence this is an element of~$I_{V_r}$. It follows that the set of $u\in V_{r-1}$ satisfying (a) is invariant under translation by $\BF_q v_0$. By the same argument with $v_0,\ldots,v_{r-1}$ interchanged it is invariant under translation by $\BF_q v_i$ for all $0\le i\le r-1$. It is therefore invariant under translation by all of~$V_{r-1}$. Since it holds for $u=0$, it follows for all $u\in V_{r-1}$, as desired.
To prove (b) observe that
$$\sum_{u\in\circV_{r-1}} \biggl( \frac{1}{X_r+u} - \frac{1}{X_r} \biggr)
\ =\ \biggl(\sum_{u\in  V_{r-1}} \frac{1}{X_r+u}\biggr) - q^{r-1}\cdot \frac{1}{X_r}
\ =\ f_r,$$
because $q^{r-1}=0$ in~$\BF_q$. Thus summing (a) over all $u\in\circV_{r-1}$ implies that $\frac{1}{v_0\cdots v_{r-1}} \cdot f_r \in I_{V_r}$ and hence~(b).
\end{Proof}


\begin{Lem}
\label{dual:IrEqual}
$I_{V_{r-1}}\cdot J_r = I_{V_r}$ for any $r\ge2$. 
\end{Lem}

\begin{Proof}
By Lemma \ref{dual:IrSubset} it remains to prove the inclusion `$\supset$'. Since the left hand side is an ideal of $R_{V_r}$ by Lemma \ref{dual:JrIdeal}, it suffices to show that each generator $\frac{1}{v_0\cdots v_r}$ of $I_{V_r}$ lies in $I_{V_{r-1}}\cdot J_r$. For this let $n$ denote the number of indices $0\le i\le r$ with $v_i\not\in V_{r-1}$. As any $r$ of the vectors $v_0,\ldots,v_r$ are linearly independent and $V_{r-1}$ has dimension~$r-1$, we have $2\le n\le r+1$. After renumbering we may therefore assume that $v_0$, $v_1\not\in V_{r-1}$. After multiplying these elements by suitable constants in $\BF_q^\times$ we may further assume that $v_0=X_r+u_0$ and $v_1=X_r+u_1$ for $u_0$, $u_1\in V_{r-1}$. Since $r\ge2$, the subset $\{v_0,v_1\} \subset \{v_0,\ldots,v_{r-1}\}$ is linearly independent; hence $v_0\not= v_1$. The identity (\ref{Relation3}) thus implies that
$$\frac{1}{v_0\cdots v_r}
\ =\ \frac{1}{v_1-v_0} \cdot \biggl(\frac{1}{v_0}-\frac{1}{v_1}\biggr) \cdot \frac{1}{v_2\cdots v_r}$$
with $v_1-v_0 = u_1-u_0 \in \circV_{r-1}$. If $n=2$, this element is equal to
$$\frac{1}{(v_1{-}v_0)v_2\cdots v_r} \cdot \biggl[ \biggl(\frac{1}{X_r+u_0}-\frac{1}{X_r}\biggr) 
                                             - \biggl(\frac{1}{X_r+u_1}-\frac{1}{X_r}\biggr) \biggr].$$
Here all of the $r$ factors in the first denominator lie in $V_{r-1}$ and any $r-1$ of them are linearly independent; hence the first factor lies in~$I_{V_{r-1}}$. As the second factor is a sum of generators of~$J_r$, the product lies in $I_{V_{r-1}}\cdot J_r$, as desired. If $n\ge3$, the element in question is 
$$\frac{1}{v_0(v_1{-}v_0)v_2\cdots v_r} 
- \frac{1}{v_1(v_1{-}v_0)v_2\cdots v_r}.$$
Here again any $r$ factors in each denominator are linearly independent, but the number of factors not in $V_{r-1}$ is now $n-1$. Thus the desired assertion follows by induction on~$n$.
\end{Proof}


\begin{Prop}
\label{dual:IrBasis}
For any $r\ge1$ the ideal $I_{V_r}$ is a free module of rank $q^{\frac{r(r-1)}{2}}$ with basis $\hat\Delta_1 \cdots \hat\Delta_r$ over $\BF_q[f_1, \ldots, f_r]$.
\end{Prop}

\begin{Proof}
We first claim that $I_{V_r}$ is generated as an $\BF_q[f_1, \ldots, f_r]$-module by $\hat\Delta_1\cdots \hat\Delta_r$. We will prove this using induction on~$r$. For $r=1$ the claim follows from the fact that $f_1 = \frac{1}{X_1}$ and hence $I_{V_1} = \bigl(\frac{1}{X_1^2}\bigr) = \bigl(\frac{f_r}{X_1}\bigr)$. For $r\ge2$ we have 
$$I_{V_r} \ =\ I_{V_{r-1}}\cdot J_r \ =\ \sum_{\epsilon\in \hat\Delta_r} I_{V_{r-1}}[f_r]\cdot\epsilon$$
by Lemma \ref{dual:IrEqual} and the definition of~$J_r$. By the claim for $r-1$ this is equal to
$$\sum_{\epsilon\in \hat\Delta_r} \sum_{\epsilon'\in \hat\Delta_1\cdots \hat\Delta_{r-1}} 
  \BF_q[f_1, \ldots, f_r] \cdot\epsilon\cdot\epsilon',$$
proving the claim for~$r$, as desired.

Next observe that $|\hat\Delta_r| = q^{r-1}$ for all $r\ge1$, and hence $|\hat\Delta_1\cdots \hat\Delta_r| \le q^{\frac{r(r-1)}{2}}$ for all $r\ge1$. On the other hand recall from Theorem \ref{RV2:free2} that $R_{V_r}$ is a free $\BF_q[f_1, \ldots, f_r]$-module of rank $q^{\frac{r(r-1)}{2}}$. As $I_{V_r}$ is a non-zero ideal of the integral domain $R_{V_r}$, it contains a submodule isomorphic to $R_{V_r}$, and so $I_{V_r}$ is an $\BF_q[f_1, \ldots, f_r]$-module of rank equal to $q^{\frac{r(r-1)}{2}}$. Since this is already an upper bound for the number of generators $|\hat\Delta_1\cdots \hat\Delta_r|$, this number is equal to $q^{\frac{r(r-1)}{2}}$ and the generators linearly independent over $\BF_q[f_1, \ldots, f_r]$, as desired.
\end{Proof}


Next recall from (\ref{RV2:DeltarDef}) that 
$$\Delta_r\ =\ \{1\} \cup \biggl\{ \frac{1}{X_r + u} \;\biggm|\; u \in \circV_{r-1} \biggr\}.$$
Consider the bijection $\Delta_r\to \hat\Delta_r$, $\delta\mapsto\hat\delta$ defined by
$$1\mapsto\frac{f_r}{X_r} \qquad\mbox{and}\qquad
\frac{1}{X_r + u} \mapsto \frac{1}{X_r + u} - \frac{1}{X_r}.$$
For any $r\ge2$ let $W_r$ denote the group of automorphisms of $V_r$ that fix every element of $V_{r-1}$ and send $X_r$ to $X_r+u$ for some $u\in V_{r-1}$. (Compare the proof of Proposition~\ref{RV2:free1}.) By functoriality this group acts on~$R_{V_r}$. Consider the operator
$$M_r := \sum_{\sigma\in W_r}\sigma :\ R_{V_r}\longto R_{V_r}.$$

\begin{Lem}
\label{dual:MrCalc}
For any $r\ge2$ and any $\delta$, $\delta'\in\Delta_r$ we have
$$M_r(\delta\cdot\hat\delta')\ =\ 
\left\{\begin{array}{ll}
f_r^2 &\mbox{if\ }\delta=\delta',\\
0     &\mbox{otherwise.}\\
\end{array}\right.$$
\end{Lem}

\begin{Proof}
This amounts to the equations
\begin{itemize}
\item[(a)] $\displaystyle M_r\biggl(\frac{f_r}{X_r}\biggr) = f_r^2$,
\item[(b)] $\displaystyle M_r\biggl(\frac{f_r}{(X_r+u)X_r}\biggr) = 0$ for all $u\in\circV_{r-1}$, 
\item[(c)] $\displaystyle M_r\biggl(\frac{1}{X_r+u'}-\frac{1}{X_r}\biggr) = 0$ for all $u'\in\circV_{r-1}$,
\item[(d)] $\displaystyle M_r\biggl(\frac{1}{X_r+u}\cdot\biggl(\frac{1}{X_r+u'}-\frac{1}{X_r}\biggr)\biggr) = 0$ for all $u$, $u'\in\circV_{r-1}$ with $u\not=u'$, and
\item[(e)] $\displaystyle M_r\biggl(\frac{1}{X_r+u}\cdot\biggl(\frac{1}{X_r+u }-\frac{1}{X_r}\biggr)\biggr) = f_r^2$ for all $u\in\circV_{r-1}$.
\end{itemize}
To prove them, note first that 
$$f_r\ =\ \sum_{u\in  V_{r-1}} \frac{1}{X_r+u}
\ =\ M_r\biggl(\frac{1}{X_r}\biggr) 
\ =\ M_r\biggl(\frac{1}{X_r+u'}\biggr)$$
for all $u'\in\circV_{r-1}$. This immediately yields (c) and, since $f_r$ is invariant under~$W_r$, it implies~(a). Next the identity (\ref{Relation3}) implies that
\begin{eqnarray*}
M_r\biggl(\frac{f_r}{(X_r+u)X_r}\biggr)
&=& M_r\biggl(\frac{f_r}{u}\cdot\biggl(\frac{1}{X_r}-\frac{1}{X_r+u}\biggr)\biggr) \\
&=& \frac{f_r}{u}\cdot \biggl[M_r\biggl(\frac{1}{X_r}\biggr)-M_r\biggl(\frac{1}{X_r+u}\biggr)\biggr] \\
&=& 0,
\end{eqnarray*}
whence~(b). Similarly
\begin{eqnarray*}
M_r\biggl(
\rlap{$\displaystyle 
\frac{1}{X_r+u}\cdot \biggl(\frac{1}{X_r+u'}-\frac{1}{X_r}\biggr)\biggr)$} && \\
&\!\!\stackrel{(\ref{Relation3})}{=}\!\!
& M_r\biggl( \frac{1}{u-u'}\cdot\biggl(\frac{1}{X_r+u'}-\frac{1}{X_r+u}\biggr)
             -\frac{1}{u   }\cdot\biggl(\frac{1}{X_r   }-\frac{1}{X_r+u}\biggr) \biggr) \\
&=& \frac{1}{u-u'}\cdot\biggl[M_r\biggl(\frac{1}{X_r+u'}\biggr)-M_r\biggl(\frac{1}{X_r+u}\biggr)\biggr]
   -\frac{1}{u   }\cdot\biggl[M_r\biggl(\frac{1}{X_r   }\biggr)-M_r\biggl(\frac{1}{X_r+u}\biggr)\biggr] \\
&=& 0,
\end{eqnarray*}
proving~(d). Finally, 
$$M_r\biggl(\frac{1}{X_r+u}\cdot\biggl(\frac{1}{X_r+u }-\frac{1}{X_r}\biggr)\biggr)
\ \stackrel{(\ref{Relation3})}{=}\ M_r\biggl(\frac{1}{(X_r+u)^2}\biggr)
    -M_r\biggl(\frac{1}{u}\cdot\biggl(\frac{1}{X_r}-\frac{1}{X_r+u}\biggr)\biggr).$$
Here the second summand vanishes as before, and the first summand is by the definition of~$M_r$ equal to
$$\sum_{u'\in V_{r-1}} \frac{1}{(X_r+u+u')^2} \ =\ \sum_{w\in V_{r-1}} \frac{1}{(X_r+w)^2}.$$
On the other hand the definition of $f_r$ implies that 
\begin{eqnarray*}
f_r^2 &=& \sum_{w,w'\in V_{r-1}} \frac{1}{(X_r+w)(X_r+w')} \\
&\!\!\stackrel{(\ref{Relation3})}{=}\!\!
& \sum_{w\in V_{r-1}} \frac{1}{(X_r+w)^2}
  + \sum_{w\not=w'}\frac{1}{w'-w}\cdot\biggl(\frac{1}{X_r+w}-\frac{1}{X_r+w'}\biggr) \\
&=& \sum_{w\in V_{r-1}} \frac{1}{(X_r+w)^2}
  + \sum_{w\not=w'}\frac{1}{w'-w}\cdot\frac{1}{X_r+w }
  - \sum_{w\not=w'}\frac{1}{w'-w}\cdot\frac{1}{X_r+w'}.
\end{eqnarray*}
By rearranging the last two sums we find that they cancel each other out. Together this shows~(e), finishing the proof.
\end{Proof}


For any $r\ge1$ consider the map
$$\Delta_1\cdots\Delta_r \longto \hat\Delta_1\cdots \hat\Delta_r,\ \ 
\delta=\delta_1\cdots\delta_r \longmapsto \hat\delta := \hat\delta_1\cdots\hat\delta_r,$$
which is a well-defined bijection following Theorem \ref{RV2:free2} (b) and Proposition~\ref{dual:IrBasis}. Consider the operator
$$N_r := \sum_{\sigma\in U_r}\sigma :\ R_{V_r}\longto R_{V_r},$$
where $U_r$ is the group of upper-triangular matrices defined in Theorem~\ref{RV2:invariants}.

\begin{Lem}
\label{dual:NrCalc}
For any $r\ge1$ and any $\delta$, $\delta'\in\Delta_1\cdots\Delta_r$ we have
$$N_r(\delta\cdot\hat\delta')\ =\ 
\left\{\begin{array}{ll}
f_1^2\cdots f_r^2 &\mbox{if\ }\delta=\delta',\\
0                &\mbox{otherwise.}\\
\end{array}\right.$$
\end{Lem}

\begin{Proof}
For $r=1$ we have $N_1=\id$ and $\delta=\delta'=1$ and $\hat\delta'=\frac{f_1}{X_1}=f_1^2$, whence the assertion. Suppose that $r\ge2$ and we have proved the assertion for $r-1$. Then viewing $U_{r-1}$ as a group of automorphisms of $V_r$ that fixes $X_r$ we obtain a semi-direct product decomposition $U_r = U_{r-1}\ltimes W_r$. This decomposition implies that $N_r = N_{r-1}\circ M_r$. Write $\delta=\delta_1\cdots\delta_r$ and $\delta'=\delta'_1\cdots\delta'_r$ with $\delta_i$, $\delta'_i\in\Delta_i$. Since $M_r$ is $R_{V_{r-1}}$-linear, we find that
\begin{eqnarray*}
N_r(\delta\cdot\hat\delta')
&=& N_{r-1}\circ M_r \bigl(\delta_1\cdots\delta_{r-1}\hat\delta'_1\cdots\hat\delta'_{r-1}
                           \cdot\delta_r\hat\delta'_r\bigr) \\
&=& N_{r-1}\bigl(\delta_1\cdots\delta_{r-1}\hat\delta'_1\cdots\hat\delta'_{r-1}
                 \cdot M_r(\delta_r\hat\delta'_r)\bigr).
\end{eqnarray*}
By Lemma \ref{dual:MrCalc} this is
$$= \left\{\begin{array}{ll}
N_{r-1}\bigl(\delta_1\cdots\delta_{r-1}\hat\delta'_1\cdots\hat\delta'_{r-1}\cdot f_r^2\bigr)
&\mbox{if\ }\delta_r=\delta'_r,\\
0     &\mbox{otherwise.}\\
\end{array}\right.$$
But $f_r^2$ is invariant under $U_{r-1}$; hence in the first case we get
$$N_{r-1}\bigl(\delta_1\cdots\delta_{r-1}\hat\delta'_1\cdots\hat\delta'_{r-1}\bigr)\cdot f_r^2.$$
Together the lemma follows by induction on~$r$.
\end{Proof}


\begin{Prop}
\label{dual:IrIsom}
For any $r\ge1$ there is an isomorphism of graded $R_{V_r}$-modules
$$I_{V_r} \stackrel{\sim}{\longto}
\Hom_{R_{V_r}^{U_r}} \bigl( R_{V_r} , f_1^2\cdots f_r^2 R_{V_r}^{U_r} \bigr),$$
where $R_{V_r}$ acts on the right hand side by the formula $(x\phi)(y) := \phi(xy)$.
\end{Prop}

\begin{Proof}
By definition the operator $N_r$ sends $R_{V_r}$ to the subring of $U_r$-invariants $R_{V_r}^{U_r}$ and is $R_{V_r}^{U_r}$-linear. Recall from Theorem \ref{RV2:invariants} that this subring is $\BF_q[f_1, \ldots, f_r]$. Thus Lemma \ref{dual:NrCalc} implies that we have a perfect graded $R_{V_r}^{U_r}$-bilinear pairing
$$I_{V_r}\times R_{V_r} \longto f_1^2\cdots f_r^2 R_{V_r}^{U_r},
\ (x,y) \mapsto N_r(xy).$$
It follows that the map $x\mapsto \bigl(y\mapsto N_r(xy)\bigr)$ has the desired properties.
\end{Proof}


\begin{Prop}
\label{dual:DualCrit}
Let $S$ be a Cohen-Macaulay graded algebra of finite type over a field~$k$. Let $N$ be a finitely generated graded $S$-module whose associated coherent sheaf $\tilde N$ on $\Proj S$ is a dualizing sheaf. Let $R$ be a graded $S$-algebra which is free of finite rank as an $S$-module. Then $R$ is Cohen-Macaulay and the coherent sheaf $\tilde M$ associated to the graded $R$-module $M := \Hom_S(R,N)$ with the action $(r\phi)(y) := \phi(ry)$ for $r\in R$ is a dualizing sheaf on $\Proj R$.
\end{Prop}

\begin{Proof}
The Cohen-Macaulay property for $R$ is shown in essentially the same way as in the proof of \ref{RV1:CMnormal}. Next let $\pi:X\to Y$ denote the finite morphism of schemes $\Proj R \to \Proj S$. Let $\pi^{!}(\tilde N)$ denote the sheaf of $\CO_X$-modules corresponding to the sheaf of $\pi_* \CO_X$-modules $\CHom_{\CO_Y}(\pi_* \CO_X, \tilde N)$. Then \cite[Ch.$\,$III, Excs.$\,$6.10a, 7.2a]{Hartshorne} imply that $\pi^{!}(\tilde N)$ is a dualizing sheaf for $X$, and $\pi^{!}(\tilde N) = \tilde M$ by the definition of~$M$.
\end{Proof}


\medskip
\begin{Proofof}{Theorem \ref{dual:DualQV}}
Set $S:=R_{V_r}^{U_r}$ and $R:=R_{V_r}$ and $N:=f_1\cdots f_r R_{V_r}^{U_r}$. Then $S=\BF_q[f_1, \ldots, f_r]$ and $R$ is free of finite type as $S$-module by Theorems \ref{RV2:invariants} and~\ref{RV2:free2}. In particular $S$ is Cohen-Macaulay and $\Proj S = \BP^{r-1}_{\BF_q}$. 
Also $\tilde N$ is the ideal sheaf of the equation $f_1\cdots f_r=0$, which defines a union of $r$ hyperplanes; hence $\tilde N \cong \CO(-r) \cong \omega_{\BP^{r-1}}$ is a dualizing sheaf on $\Proj S$. Thus all the assumptions in Proposition \ref{dual:DualCrit} are satisfied. 
On the other hand Proposition \ref{dual:IrIsom} implies that
\begin{eqnarray*}
M &:=& \Hom_S(R,N) \ =\ 
\Hom_{R_{V_r}^{U_r}} \bigl( R_{V_r} , 
  {\textstyle\frac{1}{f_1\cdots f_r}}f_1^2\cdots f_r^2 R_{V_r}^{U_r} \bigr) \\
&\cong& \frac{1}{f_1\cdots f_r}\cdot I_{V_r}
\ \cong\ \Hom_R(f_1\cdots f_rR,I_{V_r})
\end{eqnarray*}
as a graded $R$-module. Let $\CI_{V_r}$ be the ideal sheaf on $\Proj R = Q_{V_r}$ associated to the homogeneous ideal $I_{V_r}\subset R$. Then the coherent sheaf associated to $M$ is $\tilde M \cong \CHom(\CO_{Q_V}(-r),\CI_{V_r}) \allowbreak \cong \nobreak \CI_{V_r}(r)$. By Proposition \ref{dual:DualCrit} this is a dualizing sheaf on $Q_{V_r}$, as desired.
\end{Proofof}



\section{Modular interpretation}
\label{modonly}

Let $V$ be a non-zero finite dimensional $\BF_q$-vector space. Let $S$ be a scheme over~$\BF_q$, let $\CL$ be an invertible sheaf on $S$, and let $\Gamma(S,\CL)$ denote the space of global sections of~$\CL$. The product of elements $\ell_1,\ldots,\ell_n\in\Gamma(S,\CL)$ is an element $\ell_1\cdots\ell_n\in\Gamma(S, \CL^n)$, and the inverse of a nowhere vanishing element $\ell\in\Gamma(S,\CL)$ is an element $\ell^{-1}\in\Gamma(S,\CL^{-1})$.


\begin{Def}
\label{mod:linear}
By a \emph{linear map} $\lambda: V\to\Gamma(S,\CL)$ we mean any $\BF_q$-linear map. The set of these is denoted $\Hom(V,\CL)$.
\end{Def}

\begin{Def}
\label{mod:reciprocal}
By a \emph{reciprocal map} $\rho: \circV\to\Gamma(S, \CL)$ we mean any map satisfying
\begin{itemize}
\item[(a)] $\rho(\alpha v) = \alpha^{-1} \rho(v)$ for all $v \in \circV$ and $\alpha \in \BF_q^{\times}$, and
\item[(b)] $\rho(v) \cdot \rho(v') = \rho(v+v') \cdot (\rho(v) + \rho(v'))$ in $\Gamma(S, \CL^2)$ for all $v$, $v'\in \circV$ such that $v+v'\in\circV$.
\end{itemize}
The set of these is denoted $\Rec(\circV,\CL)$.
\end{Def}


\begin{Def}
\label{mod:fiberwise0}
Consider a set $X$ and a map $\phi\!: X\to\Gamma(S,\CL)$. If for all $s \in S$, the composite map $X \xrightarrow{\phi} \Gamma(S,\CL) \rightarrow \CL \otimes_{\CO_S}k(s)$ is 
\begin{itemize}
\item[(a)] non-zero, we call $\phi$ \emph{fiberwise non-zero};
\item[(b)] injective, we call $\phi$ \emph{fiberwise injective};
\item[(c)] non-zero at every $x\in X$, we call $\phi$ \emph{fiberwise invertible}.
\end{itemize}
\end{Def}

Thus $\phi$ is fiberwise invertible if and only if it sends all $x\in X$ to nowhere vanishing sections. A linear map $V\to\Gamma(S,\CL)$ is fiberwise injective if and only if its restriction to $\circV$ is fiberwise invertible. Unravelling Definition \ref{mod:reciprocal} we deduce:

%
%

\begin{Prop}
\label{mod:rec-lin-1}
Any fiberwise injective linear map $\lambda: V\to\Gamma(S,\CL)$ corresponds to a fiberwise invertible reciprocal map $\rho: \circV\to\Gamma(S, \CL^{-1})$ by the formula $\rho(v)=\lambda(v)^{-1}$ for all $v\in\circV$, and vice versa.
\end{Prop}


Next let $i\!:V'\into V$ be the inclusion of a non-zero $\BF_q$-subspace and $\pi\!:V\onto V''$ the projection to a non-zero $\BF_q$-quotient space of~$V$.

\begin{Def}
\label{mod:pullback}
The \emph{pullback under $\pi$} of a linear map $\lambda: V''\to\Gamma(S,\CL)$ is the linear map 
$$\pi^*\lambda:\ V\to\Gamma(S,\CL),\ \ 
v\mapsto \rlap{\mbox{$\lambda(\pi(v)).$}}
\hphantom{\begin{cases} \rho(v) & \text{if $v \in V'$,} \\ 0 & \text{otherwise.} \end{cases}}$$
\end{Def}

\begin{Def}
\label{mod:extension}
The \emph{extension by zero} of a reciprocal map $\rho: \circVprime\to\Gamma(S,\CL)$ is the map 
$$i_*\rho\!:\ \circV\to\Gamma(S,\CL), \ \ 
v \mapsto \begin{cases} \rho(v) & \text{if $v \in V'$,} \\ 0 & \text{otherwise,} \end{cases}$$
which by direct calculation is again reciprocal.
\end{Def}


\begin{Prop}
\label{mod:rec-lin-2}
Let $S$ be the spectrum of a field. (In this case it seems baroque to speak of `fiberwise' non-zero, injective, or invertible, so we drop the adverb.) 
\begin{itemize}
\item[(a)] Any non-zero linear map $V\to\Gamma(S,\CL)$ is equal to $\pi^*\lambda$ for a unique non-zero quotient $\pi\!:V\onto V''$ and a unique injective linear map $\lambda\!: V''\to\Gamma(S,\CL)$.
\item[(b)] Any non-zero reciprocal map $\circV\to\Gamma(S,\CL)$ is equal to $i_*\rho$ for a unique non-zero subspace $i\!:V'\into V$ and a unique invertible reciprocal map $\rho\!:V'\to\Gamma(S,\CL)$.
\end{itemize}
\end{Prop}

\begin{Proof}
(a) is obvious and included only for comparison. (b) is equivalent to saying that for any reciprocal map $\rho\!:\circV\to\Gamma(S,\CL)$, the set $V' := \{ 0 \} \cup \{ v \in \circV \mid \rho(v) \neq 0 \}$ is an $\BF_q$-subspace of~$V$. But Definition \ref{mod:reciprocal} (a) implies that $\BF_q^\times\cdot V'\subset V'$,
and \ref{mod:reciprocal} (b) implies that for all $v$, $v'\in\circV$ with $v+v'\in\circV$ and $\rho(v)$, $\rho(v')\not=0$ we have $\rho(v+v')\not=0$. Therefore $V'+V'\subset V'$, as desired.
\end{Proof}


Now we turn to the modular interpretation. Two pairs consisting of an invertible sheaf and a linear map $(\CL, \lambda)$ and $(\CL', \lambda')$ are called \emph{isomorphic} if there exists an isomorphism of invertible sheaves $\CL \cong \CL'$ that is compatible with $\lambda$ and~$\lambda'$. Similarly for reciprocal maps.

\medskip
Recall from Section \ref{RV1} that $P_V=\Proj S_V$ where $S_V$ is the symmetric algebra of $V$ over~$\BF_q$. The natural isomorphism $\lambda_V\!:V = S_{V,1} \xrightarrow{\sim} \Gamma(P_V,\CO_{P_V}(1))$ is then a fiberwise non-zero linear map. The well-known description of the functor of points of projective space \cite[Ch.$\,$II Thm.$\,$7.1]{Hartshorne} yields:

\begin{Prop}
\label{mod:PV}
The scheme $P_V$ with the universal family $(\CO_{P_V}(1),\lambda_V)$ represents the functor which associates to a scheme $S$ over $\BF_q$ the set of isomorphism classes of pairs $(\CL, \lambda)$ consisting of an invertible sheaf $\CL$ on $S$ and a fiberwise non-zero linear map $\lambda\!:V\to\Gamma(S,\CL)$.
\end{Prop}

The description of $\Omega_V$ from Section \ref{RV1} implies:

\begin{Prop}
\label{mod:omega1}
The open subscheme $\Omega_V\subset P_V$ represents the subfunctor of fiberwise injective linear maps.
\end{Prop}

On the other hand consider the natural map $\rho_V\!: \circV \to R_{V,-1} \cong \Gamma(Q_V, \CO_{Q_V}(1))$ given by $v \mapsto \frac{1}{v}$. The identities (\ref{Relation1}) and (\ref{Relation2}) show that $\rho_V$ is reciprocal in the sense of Definition~\ref{mod:reciprocal}. It is also fiberwise non-zero, because the elements $\frac{1}{v}$ generate the augmentation ideal of~$R_V$.

\begin{Thm}
\label{mod:QV}
The scheme $Q_V$ with the universal family $(\CO_{Q_V}(1),\rho_V)$ represents the functor which associates to a scheme $S$ over $\BF_q$ the set of isomorphism classes of pairs $(\CL, \rho)$ consisting of an invertible sheaf $\CL$ on $S$ and a fiberwise non-zero reciprocal map $\rho\!:\circV\to\Gamma(S,\CL)$.
\end{Thm}

\begin{Proof}
Recall from Section \ref{RV1} that $A_V$ is the polynomial ring over $\BF_q$ in the indeterminates $Y_v$ for all $v\in\circV$. Consider the map $\sigma_V\!: \circV\to \Gamma(\Proj A_V,\CO(1))$, $v\mapsto Y_v$. Then  
\cite[Ch.$\,$II Thm.$\,$7.1]{Hartshorne} says that $\Proj A_V$ with the universal family $(\CO(1),\sigma_V)$ represents the functor of isomorphism classes of pairs $(\CL,\rho)$ consisting of an invertible sheaf $\CL$ on $S$ and a fiberwise non-zero map $\rho\!:\circV\to\Gamma(S,\CL)$. Theorem \ref{RV1:presentation} implies that $Q_V$ is the closed subscheme of $\Proj A_V$ defined by the homogenous relations~\ref{RV1:generators}. By Definition \ref{mod:reciprocal} these relations are precisely those that require $\rho$ to be reciprocal.
\end{Proof}

\medskip
Proposition \ref{mod:rec-lin-1} implies:

\begin{Prop}
\label{mod:omega2}
The open subscheme $\Omega_V\subset Q_V$ represents the subfunctor of fiberwise invertible reciprocal maps.
\end{Prop}


For use in the next section we include the following variant. Let $\tilde Q_V := \Spec R_V$ be the affine cone over~$Q_V$, and consider the reciprocal map $\tilde\rho_V\!: \circV \into R_V \cong \Gamma(\tilde Q_V, \CO_{\tilde Q_V})$ given by $v \mapsto \frac{1}{v}$. The same reasoning as in the preceding proof shows:

\begin{Thm}
\label{mod:tildeQV}
The scheme $\tilde Q_V$ with the universal reciprocal map $\tilde\rho_V$ represents the functor which associates to a scheme $S$ over $\BF_q$ the set of all reciprocal maps $\circV\to\Gamma(S,\CO_S)$.
\end{Thm}



\section{Stratification}
\label{strat}

We keep the notations of Section \ref{modonly}. Let $\pi\!:V\onto V''$ be the projection to a non-zero $\BF_q$-quotient space. Since the pullback $\pi^*$ from Definition \ref{mod:pullback}  of a fiberwise non-zero linear map is again fiberwise non-zero, it defines a morphism of functors and hence a morphism of moduli schemes $P_{V''}\to P_V$. For simplicity we denote this morphism again by~$\pi^*$. It is a closed embedding onto an $\BF_q$-rational linear subspace of~$P_V$. Consider the composite locally closed embedding 
$$\Omega_{V''} \longinto P_{V''} \stackrel{\pi^*}{\longinto} P_V.$$
Its image represents the subfunctor of linear maps on $V$ which factor through fiberwise injective linear maps on~$V''$. We identify $\Omega_{V''}$ and $P_{V''}$ with their images. 

\begin{Thm}
\label{mod:PV-strat}
\begin{itemize}
\item[(a)] The scheme $P_V$ is the set-theoretic disjoint union of the locally closed subschemes $\Omega_{V''}$ for all non-zero quotients $V''$ of~$V$.
\item[(b)] The closure of $\Omega_{V''}$ in $P_V$ is the union of the $\Omega_{W''}$ for all non-zero quotients $W''$ of~$V''$.
\end{itemize}
\end{Thm}

\begin{Proof}
Proposition \ref{mod:rec-lin-2}~(a) implies that every point on $P_V$ lies in $\Omega_{V''}$ for a unique quotient~$V''$, which shows~(a). Part (b) follows from the fact that $\Omega_{V''}$ is open and dense in the projective scheme $P_{V''}$ and by (a) applied to~$P_{V''}$.
\end{Proof}

\medskip
We now proceed analogously for~$Q_V$. Let $i\!:V'\into V$ be the inclusion of a non-zero $\BF_q$-subspace of~$V$. Since the extension by zero $i_*$ from Definition \ref{mod:extension}  of a fiberwise non-zero reciprocal map is again fiberwise non-zero, it defines a morphism of functors and hence a morphism of moduli schemes $Q_{V'}\to Q_V$. For simplicity we denote this morphism again by~$i_*$. Its image is the subfunctor of fiberwise non-zero reciprocal maps $\rho: \circV \rightarrow \Gamma(S,\CL)$ satisfying $\rho(v)=0$ for all $v\in V\setminus V'$. As this is a closed condition, the morphism $i_*$ is a closed embedding. Consider the composite locally closed embedding 
$$\Omega_{V'} \longinto Q_{V'} \stackrel{i_*}{\longinto} Q_V.$$
Its image represents the subfunctor of reciprocal maps on $\circV$ which are the extension by zero of fiberwise invertible reciprocal maps on~$\circVprime$. We identify $\Omega_{V'}$ and $Q_{V'}$ with their images.

\begin{Thm}
\label{mod:QV-strat}
\begin{itemize}
\item[(a)] The scheme $Q_V$ is the set-theoretic disjoint union of the locally closed subschemes $\Omega_{V'}$ for all non-zero subspaces $V'$ of~$V$.
\item[(b)] The closure of $\Omega_{V'}$ in $Q_V$ is the union of the $\Omega_{W'}$ for all non-zero subspaces $W'$ of~$V'$.
\end{itemize}
\end{Thm}

\begin{Proof}
Proposition \ref{mod:rec-lin-2}~(b) implies that every point on $Q_V$ lies in $\Omega_{V'}$ for a unique subspace~$V'$, which shows~(a). Part (b) follows from the fact that $\Omega_{V'}$ is open and dense in the projective scheme $Q_{V'}$ and by (a) applied to~$Q_{V'}$.
\end{Proof}


\medskip
Now we determine the local structure of $Q_V$ along the stratum $\Omega_{V'}$ associated to a non-zero proper subspace~$V'$. Consider the subfunctor of the functor represented by $Q_V$ consisting of reciprocal maps whose restriction to $\circVprime$ are fiberwise invertible. This is an open condition; hence the subfunctor is represented by an open subscheme $U^{V'}_V \subset Q_V$. Moreover, the restriction of reciprocal maps to $\circVprime$ induces a morphism of functors and hence of schemes
$$i^*\!:\ U^{V'}_V \longto \Omega_{V'}.$$
Furthermore, we have $\Omega_{V'} \subset U^{V'}_V$, and the restriction of $i^*$ to $\Omega_{V'}$ is the identity.

\medskip
Next we fix a subspace $V''\subset V$ complementary to $V'$ and an element $v'_0\in\circVprime$. Then for any reciprocal map $\rho\!: \circV\to\Gamma(S,\CL)$ whose restriction to $\circVprime$ is fiberwise invertible, the section $\rho(v'_0)$ vanishes nowhere and thus induces an isomorphism $\CO_S\stackrel{\sim}{\to}\CL$. The composite map
$$n(\rho)\!:\ \ \circVpprime \longto \Gamma(S,\CL) \stackrel{\sim}{\longto} \Gamma(S,\CO_S),\ \ 
v'' \mapsto \frac{\rho(v'')}{\rho(v'_0)}$$
is then invariant under isomorphisms of the pair $(\CL,\rho)$. It is clearly reciprocal; hence in view of Theorem \ref{mod:tildeQV} the map $(\CL,\rho)\mapsto n(\rho)$ defines a morphism of functors and hence a morphism of moduli schemes 
$$n\!:\ U^{V'}_V \longto \tilde Q_{V''}.$$
Let $0\in \tilde Q_{V''}$ correspond to the identically zero reciprocal map $\circVpprime \to\Gamma(S,\CO_S)$. Then we have the following commutative diagram:
$$\xymatrix@C+20pt{
U^{V'}_V  \ar[r]^-{i^*\times n}    
& \Omega_{V'} \times \tilde Q_{V''} \\
\Omega_{V'} \ar[r]^-{\id\times\{0\}}_-{\sim} \ar@{^{ (}->}[u]^{}
& \Omega_{V'} \times \{0\} \ar@{^{ (}->}[u]^{} \\}$$

\begin{Prop}
\label{mod:localgeometry}
The morphism $i^*\times n$ induces an isomorphism from some neighborhood of $\Omega_{V'}$ in $U^{V'}_V$ to some neighborhood of $\Omega_{V'} \times \{0\}$ in $\Omega_{V'} \times \tilde Q_{V''}$.
\end{Prop}

\begin{Proof}
Consider the subfunctor of the functor represented by $U^{V'}_V$ over which the section $\rho(v')+\rho(v'')$ vanishes nowhere for all $v'\in\circVprime$ and $v''\in\circVpprime$. This is an open condition; hence the subfunctor is represented by an open subscheme $U_1 \subset U^{V'}_V$. Since $\rho(v')$ already vanishes nowhere and $\rho(v'')$ is everywhere zero on $\Omega_{V'}$, we also have $\Omega_{V'}\subset U_1$.

\medskip
On the other hand consider a fiberwise invertible reciprocal map $\rho'\!: \circVprime\to\Gamma(S,\CL)$ in the functor represented by $\Omega_{V'}$ and a reciprocal map $\rho''\!: \circVpprime\to\Gamma(S,\CO_S)$ in the functor represented by~$\tilde Q_{V''}$. On $\Omega_{V'} \times \tilde Q_{V''}$ consider the subfunctor over which $\rho'(v')+\rho'(v'_0)\cdot\rho''(v'')$ vanishes nowhere for all $v'\in\circVprime$ and $v''\in\circVpprime$. This is an open condition; hence the subfunctor is represented by an open subscheme $U_2 \subset \Omega_{V'} \times \tilde Q_{V''}$. Since $\rho'(v')$ already vanishes nowhere on $\Omega_{V'}$  and $\rho''(v'')=0$ over $\{0\}$, we also have $\Omega_{V'}\times\{0\}\subset U_2$.

\medskip
We claim that $i^*\times n$ induces an isomorphism $U_1\to U_2$. 

\medskip
Indeed, for any reciprocal map $\rho\!: \circV\to\Gamma(S,\CL)$ whose restriction to $\circVprime$ is fiberwise invertible, set $\rho' := i^*\rho\!: \circVprime\to\Gamma(S,\CL)$ and $\rho'' := n(\rho)\!: \circVpprime\to\Gamma(S,\CO_S)$. Then for all $v'\in\circVprime$ and $v''\in\circVpprime$ we have
$$\rho(v')+\rho(v'') \ =\ \rho'(v')+\rho'(v'_0)\cdot\rho''(v'').$$
The defining conditions for both $U_1$ and $U_2$ require precisely that these sections be fiberwise non-zero. This implies firstly that $i^*\times n$ sends $U_1$ to~$U_2$. Secondly, using \ref{mod:reciprocal} (b) it shows that
$$\rho(v'+v'') 
\ =\ \frac{\rho (v') \cdot \rho(v'')}                   {\rho(v') + \rho(v'')}
\ =\ \frac{\rho'(v') \cdot \rho'(v'_0)\cdot \rho''(v'')}{\rho'(v')+\rho'(v'_0)\cdot\rho''(v'')}.$$
Together with the equalities $\rho(v')=\rho'(v')$ and $\rho(v'')=\rho'(v'_0)\rho''(v'')$ this recovers $\rho$ completely from $\rho'$ and~$\rho''$. Conversely, these formulas associate to any pair of reciprocal maps $\rho': \circVprime\to\Gamma(S,\CL)$ and $\rho''\!: \circVpprime\to\Gamma(S,\CO_S)$ satisfying the condition for $U_2$ a map $\rho\!: \circV\to\Gamma(S,\CL)$. We leave it to the careful reader to verify that this map is reciprocal and satisfies the condition for~$U_1$. This then finishes the proof.
\end{Proof}


\begin{Thm}
\label{mod:singularlocus}
The singular locus of $Q_V$ is the union of all strata $\Omega_{V'}$ of codimension $\ge2$, that is, with $\dim(V/V') \geq 2$.
\end{Thm}

\begin{Proof}
The open stratum $\Omega_V$ is smooth, so consider the stratum $\Omega_{V'}$ associated to a non-zero proper subspace $V'\subset V$. Proposition \ref{mod:localgeometry} implies that a point $p \in \Omega_{V'}$ is regular in $Q_V$ if and only if $p\times\{0\}$ is regular in $\Omega_{V'}\times\tilde Q_{V''}$. Since $\Omega_{V'}$ is smooth, this is equivalent to the vertex $\{0\}$ being a regular point of~$\tilde Q_{V''}$. 
But the local ring of $\tilde Q_{V''}$ at $0$ is the localization of $R_{V''}$ at the augmentation ideal, and its associated graded ring is therefore again isomorphic to~$R_{V''}$. Thus $p$ is a regular point if and only if $R_{V''}$ is isomorphic to a polynomial ring. By Remark \ref{RV1:UFD} that is the case if and only if $\dim V''=1$.
\end{Proof}

\begin{Prop}
\label{mod:divisor}
The divisor of the section $\frac{1}{v} \in R_{V,-1} = H^0(Q_V,\CO_{Q_V}(1))$ for any $v\in\circV$ is the sum of $Q_{V'}$ for all $V'\subset V$ of codimension~$1$ with $v\not\in V'$, with multiplicity $1$ each.
\end{Prop}

\begin{Proof}
Since $\frac{1}{v}$ is invertible over the open stratum~$\Omega_V$, its divisor is a linear combination of the irreducible components of $Q_V\setminus\Omega_V$. Theorem \ref{mod:QV-strat} implies that these irreducible components are precisely the $Q_{V'}$ for all $V'$ of codimension~$1$. We fix such a $V'$ and determine the multiplicity of~$Q_{V'}$. In the case $v\in V'$ the section $\frac{1}{v}$ remains invertible over~$\Omega_{V'}$; hence the multiplicity is~$0$. Otherwise we can apply Proposition \ref{mod:localgeometry} with $V'' := \BF_q\cdot v$, in which case $R_{V''} = \BF_q\bigl[\frac{1}{v}\bigr]$ and the multiplicity is therefore~$1$.
\end{Proof}

\begin{Prop}
\label{mod:ideal}
The ideal sheaf $\CI_V$ from Section \ref{dual} coincides with $\CO_{Q_V}$ over the open stratum $\Omega_V$ and has multiplicity $2$ along all strata $\Omega_{V'}$ of codimension~$1$.
\end{Prop}

\begin{Proof}
Set $r:=\dim V$ and recall that $\CI_V$ is the ideal sheaf associated to the homogeneous ideal $I_V \subset R_V$ that is generated by the elements $\frac{1}{v_0\cdots v_r}$ for all $v_0,\ldots,v_r\in \circV$, any $r$ of which are linearly independent. Proposition \ref{mod:divisor} implies that the divisor of each generator is
$$\mathop{\rm div}\bigl({\textstyle\frac{1}{v_0\cdots v_r}}\bigr)
\ =\ \sum_{\dim(V/V')=1} \bigl|\{0\le i\le r\;|\;v_i\not\in V'\}\bigr| \cdot Q_{V'}.$$
Thus we need to show that for any fixed subspace $V'$ of codimension~$1$, the minimum of the number 
$\bigl|\{0\le i\le r\;|\;v_i\not\in V'\}\bigr|$ for all $v_0,\ldots,v_r$ as above is~$2$. 
But since any $r$ vectors in $V'$ are linearly dependent, this number is at least~$2$. On the other hand, take any basis $v_1,\ldots,v_{r-1}$ of $V'$ and any $v_r\in V\setminus V'$ and set $v_0 := v_1+\ldots+v_r$. Then any $r$ of the vectors $v_0,\ldots,v_r$ are linearly independent and the number of those not in $V'$ is~$2$. Thus the minimum is indeed~$2$, as desired.
\end{Proof}

\begin{Exer}
\label{mod:exer1}
Let $j$ denote the open embedding of the regular locus $Q_V^{\rm reg} \into Q_V$. Then $\CI_V=j_*j^*\CI_V$.
\end{Exer}

\begin{Exer}
\label{mod:exer2}
The reduced closed subscheme of $Q_V$ supported on $Q_V{\setminus}\Omega_V$ is the subscheme associated to the homogeneous ideal of $R_V$ that is generated by the elements $\frac{1}{v_1\cdots v_r}$ for all linearly independent $v_1,\ldots,v_r\in \circV$, where $r=\dim V$.
\end{Exer}



\section{Strange morphisms}
\label{maps}

The morphisms defined in this section are intriguing, but not used elsewhere in this article. Let $\Frob_{q^n}$ denote the $q^n$-th power Frobenius morphism on any scheme over~$\BF_q$. 
As before we let $r:=\dim V$. We denote the natural pairing of $V$ with its dual space $V^*$ by
$$V^* \times V \rightarrow \BF_q,\ \ (\ell,v) \mapsto \langle \ell,v \rangle.$$

\begin{Prop}
\label{maps:gfcomposition}
Consider a reciprocal map $\rho\!: \circV \rightarrow \Gamma(S,\CL)$ and a linear map $\lambda: V^* \rightarrow \Gamma(S,\CL)$.
\begin{itemize}
\item[(a)] The map $\displaystyle g_V(\rho): V^* \longrightarrow \Gamma(S,\CL), \ 
\ell \longmapsto \sum_{\myatop{v \in \circV}{\langle \ell,v \rangle =1}} \rho(v)$ is linear.
\item[(b)] The map $\displaystyle f_V(\lambda): \circV \longrightarrow \Gamma\bigl(S,\CL^{q^{r-1}}\bigr), \ 
v \longmapsto \prod_{\myatop{\ell \in V^*}{\langle \ell,v \rangle =1}} \lambda(\ell)$ is reciprocal.
\item[(c)] We have $(g_V \circ f_V)(\lambda) = \lambda^{q^{r-1}}$.
\item[(d)] We have $(f_V \circ g_V)(\rho)    = \rho   ^{q^{r-1}}$.
\end{itemize}
\end{Prop}

\begin{Proof}
The condition \ref{mod:reciprocal}~(a) implies that each summand of the sum
$$\sum_{v \in \circV / \BF_q^{\times}} \langle \ell,v \rangle \cdot \rho(v)$$
depends only on $\BF_q^\times\cdot v$; hence the sum is well-defined. As all summands with $\langle \ell,v \rangle=\nobreak0$ vanish, the sum gives an equivalent formula for $g_V(\rho)$.
This formula is linear in~$\ell$, proving~(a).
The remaining assertions rely on lengthy elementary calculations which we leave to the interested reader.
\end{Proof}

\begin{Prop}
\label{maps:fiberwise}
A reciprocal map $\rho$ is fiberwise non-zero (resp.\ fiberwise invertible) if and only if $g_V(\rho)$ is fiberwise non-zero (resp.\ fiberwise injective). A linear map $\lambda$ is fiberwise non-zero (resp.\ fiberwise injective) if and only $f_V(\lambda)$ is fiberwise non-zero (resp.\ fiberwise invertible).
\end{Prop}

\begin{Proof}
If $\rho$ is identically zero in some fiber, then $g_V(\rho)$ is identically zero in the same fiber. Conversely, if $g_V(\rho)$ is identically zero in some fiber, then $(f_V \circ g_V)(\rho) = \rho^{q^{r-1}}$ and hence $\rho$ is identically zero in the same fiber. Together this proves that $\rho$ is fiberwise non-zero if and only if $g_V(\rho)$ is fiberwise non-zero. In the same way one shows that $\lambda$ is fiberwise non-zero if and only if $f_V(\lambda)$ is fiberwise non-zero. Next, the definition of $f_V$ implies that $\lambda$ is fiberwise injective if and only if $f_V(\lambda)$ is fiberwise invertible. Applying this to $\lambda=g_V(\rho)$ we deduce that $g_V(\rho)$ is fiberwise injective if and only if $(f_V \circ g_V)(\rho) = \rho^{q^{r-1}}$ is fiberwise invertible if and only if $\rho$ is fiberwise invertible.
\end{Proof}

\begin{Thm}
\label{maps:strange}
The constructions in Proposition \ref{maps:gfcomposition} induce morphisms forming a commutative diagram:
$$\xymatrix{
\Omega_V \ar[rr]^{g_V} \ar@{^{ (}->}[d]  & &   \Omega_{V^*} \ar@{^{ (}->}[d] \ar[rr]^{f_V}   & &  \Omega_V \ar@{^{ (}->}[d]   \\
Q_V \ar[rr]^{g_V}    & &      P_{V^*}   \ar[rr]^{f_V}      & &   Q_V \\
}$$
Both composites $f_V \circ g_V$ and $g_V \circ f_V$ are the Frobenius morphism $\Frob_{q^{r-1}}$. In particular, the morphisms $f_V$ and $g_V$ are finite, bijective on the underlying sets, and radicial.
\end{Thm}

\begin{Proof}
Proposition \ref{maps:gfcomposition} (a--b) and Proposition \ref{maps:fiberwise} imply that $g_V$ and $f_V$ induce morphisms of functors and thus of schemes making the diagram commute. Proposition \ref{maps:gfcomposition} (c--d) implies the remaining assertions.
\end{Proof}

\medskip

Next let $i\!:V'\into V$ be the inclusion of a non-zero $\BF_q$-subspace of $V$ and let $\pi\!:V^*\onto V^{\prime*}$ be the projection dual to~$i$. Let $\pi^*$ be the pullback of linear maps from Definition \ref{mod:pullback} and $i_*$ the extension by zero of reciprocal maps from Definition~\ref{mod:extension}. Set $r' := \dim V'$ and $r'' := \dim V - \dim V'$. By expanding all definitions involved we find:

\begin{Prop}
\label{maps:gfcompat}
The following diagram commutes:
$$\xymatrix@C+20pt{
\Rec(\circVprime,\CL) \ar[r]^-{g_{V'}}  \ar[d]_{i_*} & 
\Hom({V'}^*, \CL) \ar[d]_{\pi^*} \ar[r]^-{f_{V'}}   & 
\Rec(\circVprime,\CL^{q^{r'-1}})  \ar[d]^{i_*\circ\Frob_{q^{r''}}}     \\
\Rec(\circV,\CL) \ar[r]^-{g_V}  & 
\Hom(V^*, \CL) \ar[r]^-{f_V}  & 
\Rec(\circV,\CL^{q^{r-1}}) \\ }$$
\end{Prop}

As a direct consequence of this and Theorem \ref{maps:strange} for $V'$ in place of $V$ we obtain:

\begin{Thm}
\label{maps:commute}
The following diagram commutes:
$$\xymatrix{
\Omega_{V'} \ar[rr]^{g_{V'}} \ar@{^{ (}->}[d]  & &   \Omega_{V^{\prime*}} \ar@{^{ (}->}[d] \ar[rr]^{f_{V'}}   
& &  \Omega_{V'} \ar@{^{ (}->}[d]   \\
Q_{V'} \ar[rr]^{g_{V'}} \ar@{^{ (}->}[d]_{i_*} & & P_{{V'}^*}  \ar@{^{ (}->}[d]_{\pi^*}   \ar[rr]^{f_{V'}}  & & Q_{V'}  \ar[d]^{i_*\circ\Frob_{q^{r''}}}  \\
Q_V \ar[rr]^{g_V}   & &   P_{V^*} \ar[rr]^{f_V}    & &    Q_V
}$$
\end{Thm}

Thus the morphisms $g_V$ and $f_V$ give a precise correspondence between the stratifications of $P_{V^*}$ and of $Q_V$ described in Section \ref{strat}.



\section{Desingularization of $Q_V$}
\label{BV}

As before let $V$ be a non-zero finite dimensional $\BF_q$-vector space. All tensor products and all fiber products will be taken over~$\BF_q$. 

\medskip
Recall from Proposition \ref{mod:PV} that $P_V$ represents the functor of isomorphism classes of pairs $(\CL, \lambda)$ consisting of an invertible sheaf $\CL$ on $S$ and a fiberwise non-zero linear map $\lambda\!:V\to\Gamma(S,\CL)$. For any such pair let $\CE_V$ denote the kernel of the surjection of coherent sheaves $\lambda\otimes\id\!: V\otimes\CO_S\onto\CL$. Note that $\CE_V$ is a locally free coherent sheaf of rank $\dim(V)-1$ which is locally a direct summand of $V\otimes\CO_S$. Clearly $\CE_V$ as a subsheaf depends only on the isomorphism class of $(\CL,\lambda)$, and conversely, $\CE_V$ determines that isomorphism class because $\CL \cong (V\otimes\CO_S)/\CE_V$. Thus we find:

\begin{Prop}
\label{BV:PV1}
The scheme $P_V$ represents the functor which to a scheme $S$ over $\BF_q$ associates the set of coherent subsheaves $\CE_V\subset V\otimes\CO_S$ such that $(V\otimes\CO_S)/\CE_V$ is locally free of rank~$1$.
\end{Prop}

Let us briefly discuss how to specify open and closed conditions on a moduli scheme. Consider coherent sheaves $\CF_1$, $\CF_2\subset\CF$ on a scheme~$S$, such that $\CF_1$ is locally free and $\CF/\CF_2$ is locally free of rank~$1$. Then the composite homomorphism $\CF_1\into\CF\onto\CF/\CF_2$ can be given in local coordinates by a tuple of local sections of~$\CO_S$. The ideal generated by these local sections is independent of the local coordinates and defines the unique maximal closed subscheme $S'\subset S$ over which the homomorphism vanishes. In other words $S'$ is the unique maximal \emph{closed} subscheme of $S$ such that $\CF_1|_{S'} \subset \CF_2|_{S'}$. Moreover, the complement $S\setminus S'$ is the set of points $s\in S$ with residue field $k(s)$ such that $\CF_1\otimes k(s) \not\subset \CF_2\otimes k(s)$. Since $\CF/\CF_2$ is locally free of rank~$1$, this is equivalent to $\CF_1\otimes k(s) + \CF_2\otimes k(s) = \CF\otimes k(s)$ and thus by Nakayama's lemma to $(\CF_1+\CF_2)_s = \CF_s$. Therefore $S\setminus S'$ is the unique maximal \emph{open} subscheme of $S$ such that $(\CF_1+\CF_2)|_{S\setminus S'} = \CF|_{S\setminus S'}$. Finally, both conditions $\CF_1\subset\CF_2$ and $\CF_1+\CF_2=\CF$ are invariant under pullback. Thus if a scheme $M$ represents a functor whose data involves sheaves $\CF_1$, $\CF_2\subset\CF$ as above, the two conditions define subfunctors which, by applying the preceding arguments to the universal family, are represented by certain complementary closed, resp.\ open subschemes of~$M$.

\begin{Prop}
\label{BV:PV2}
The open subscheme $\Omega_V \subset P_V$ represents the subfunctor of all $\CE_V$ for which $V \otimes \CO_S = \CE_V + (V' \otimes \CO_S)$ for all $0 \neq V' \subset V$.
\end{Prop}

\begin{Proof}
By the preceding remarks, the condition for any fixed~$V'$ defines the complement of the closed subscheme that represents the functor of all $\CE_V$ satisfying $V' \otimes \CO_S \subset \CE_V$. This subfunctor translates into the subfunctor of all linear maps that factor through $V/V'$, which is represented by the subscheme $P_{V'/V}$. Since $\Omega_V$ is the complement of the union of all these $P_{V'/V}$, the proposition follows.
\end{Proof}

\medskip
Now we consider the cartesian product of $P_{V'}$ for all $0\not=V'\subset V$, which represents tuples $\CE_\bullet = (\CE_{V'})_{V'}$. Then there exists a unique closed subscheme 
\addtocounter{Thm}{1}
\begin{equation}
\label{BV:BVdef}
B_V \ \subset\ \prod_{0\not=V'\subset V}P_{V'}
\end{equation}
representing the subfunctor of all $\CE_\bullet$ satisfying the closed condition
\addtocounter{Thm}{1}
\begin{equation}
\label{BV:BVcond}
\strut\rlap{\mbox{$\CE_{V''}\subset\CE_{V'}$}}%
\phantom{\mbox{$V' \otimes \CO_S = \CE_{V'} + (V'' \otimes \CO_S)$}}
\quad \rlap{\hbox{for all $0\not=V''\subset V'\subset\nobreak V$.}}
\phantom{\hbox{exists no $W\in\CF$ with $V''\subset W$ and $V'\not\subset W$.}}
\end{equation}
This will be the desired desingularization of~$Q_V$. 
First we construct a natural stratification of~$B_V$. Let $\CF$ be a flag of $V$, that is, a set of subspaces of $V$ which is totally ordered by inclusion and contains $0$ and~$V$. Then there exists a unique closed subscheme $B_\CF \subset B_V$ representing the subfunctor of all $\CE_\bullet$ satisfying (\ref{BV:BVcond}) and the closed condition
\addtocounter{Thm}{1}
\begin{equation}
\label{BV:BFcond}
\strut\rlap{\mbox{$V''\otimes\CO_S \subset \CE_{V'}$}}%
\phantom{\mbox{$V' \otimes \CO_S = \CE_{V'} + (V'' \otimes \CO_S)$}}
\quad \rlap{\vtop{\hbox{for all $0\not=V''\subset V'\subset\nobreak V$ such that there}
            \hbox{exists    $W\in\CF$ with $V''\subset W$ and $V'\not\subset W$.}}}
\phantom{\hbox{exists no $W\in\CF$ with $V''\subset W$ and $V'\not\subset W$.}}
\end{equation}
On the other hand, there exists a unique open subscheme $U_{\CF} \subset B_V$ representing the subfunctor of all $\CE_\bullet$ satisfying (\ref{BV:BVcond}) and the open condition
\addtocounter{Thm}{1}
\begin{equation}
\label{BV:UFcond}
\strut\rlap{\mbox{$V' \otimes \CO_S = \CE_{V'} + (V'' \otimes \CO_S)$}}%
\phantom{\mbox{$V' \otimes \CO_S = \CE_{V'} + (V'' \otimes \CO_S)$}}
\quad \vtop{\hbox{for all $0\not=V''\subset V'\subset\nobreak V$ such that there}
            \hbox{exists no $W\in\CF$ with $V''\subset W$ and $V'\not\subset W$.}}
\end{equation}
The locally closed subscheme $\Omega_\CF := B_\CF\cap U_\CF$ will be the stratum associated to~$\CF$.


\begin{Lem}
\label{BV:flagcomp}
For any two flags $\CF$ and $\CF'$ of $V$ we have:
\begin{itemize}
\item[(a)] $B_{\CF'} \subset B_\CF$ and $U_\CF \subset U_{\CF'}$ if $\CF\subset\CF'$,
\item[(b)] $B_\CF \cap U_{\CF'} = \emptyset$ if $\CF\not\subset\CF'$.
\item[(c)] $\Omega_\CF \cap \Omega_{\CF'} = \emptyset$ if $\CF\not=\CF'$.
\end{itemize}
\end{Lem}

\begin{Proof}
(a) is a direct consequence of the definition. For (b) take any $W\in\CF\setminus\CF'$. Then $W\not=0,V$. Set $V'':= W$ and let $V'\in\CF'$ be minimal with $W\subset V'$. Then $W\subsetneqq V'$, so that the condition (\ref{BV:BFcond}) applies to the subspaces $V''\subset V'$ and the filtration~$\CF$; while the condition (\ref{BV:UFcond}) applies to the subspaces $V''\subset V'$ and the filtration~$\CF'$. Thus on $B_\CF \cap U_{\CF'}$ we simultaneously have $V''\otimes\CO_S \subset \CE_{V'}$ and $V' \otimes \CO_S = \CE_{V'} + V'' \otimes \CO_S$, and hence $\CE_{V'} = V' \otimes \CO_S$, which contradicts the assumption on $\CE_{V'}$ unless $S=\emptyset$. This proves (b). Finally, (b) implies that $\Omega_\CF \cap \Omega_{\CF'} = \emptyset$ if  $\CF\not\subset\CF'$. By symmetry this yields~(c).
\end{Proof}

\begin{Thm}
\label{BV:stratification}
As a set $B_V$ is the disjoint union of the strata $\Omega_\CF$ for all flags $\CF$ of~$V$.
\end{Thm}

\begin{Proof}
The disjointness is Lemma \ref{BV:flagcomp} (c). It remains to see that every point on $B_V$ lies in $\Omega_\CF$ for some flag~$\CF$. For this we take any tuple $\CE_\bullet$ defined over the spectrum of a field~$k$ and  satisfying (\ref{BV:BVcond}). Then $\CE_\bullet$ corresponds to a collection of $k$-subspaces $E_{V'} \subset V'_k := V'\otimes k$ of codimension~$1$ such that $E_{V''}\subset E_{V'}$ for all $0\not=V''\subset V'\subset V$. We must find a flag $\CF$ of $V$ that satisfies (\ref{BV:BFcond}) and (\ref{BV:UFcond}), i.e., such that for all $0\not=V''\subset V'\subset V$:
\addtocounter{Thm}{1}
\begin{equation}
\label{BV:OFcondk}
\begin{cases}
V''_k \subset \rlap{\hbox{$E_{V'}$}} \phantom{E_{V'} + V''_k}
& \hbox{if there exists    $W\in\CF$ with $V''\subset W$ and $V'\not\subset W$,} \\[4pt]
\rlap{\hbox{$V'_k$}}
\phantom{V''_k} = E_{V'} + V''_k & \hbox{if there exists no $W\in\CF$ with $V''\subset W$ and $V'\not\subset W$.}
\end{cases}
\end{equation}
As a preparation observe that since $E_{V'} \subset V'_k$ has codimension~$1$ we always have
$$V'_k = E_{V'} + V''_k \quad\hbox{if and only if}\quad V''_k\not\subset E_{V'}.$$
Suppose first that for all non-zero $\BF_q$-subspaces $U \subset V$ we have $U_k\not\subset E_V$. Then for all $0\not=V''\subset V'\subset V$ we have $V''_k\not\subset E_{V'}\subset E_V$ and hence the second case of~(\ref{BV:OFcondk}). Thus the trivial flag $\{0,V\}$ does the job in this case.

\medskip
Otherwise there exist $\BF_q$-subspaces $0\neq U \subset V$ with $U_k\subset E_V$. Their sum $U_{\rm max}$ then enjoys the same properties and is therefore the unique largest one among them. Since $E_V\subset V_k$ has codimension~$1$, we have $U_{\rm max}\not=V$. By induction on $\dim V$ we may therefore assume that there exists a flag $\CF'$ of $U_{\rm max}$ such that (\ref{BV:OFcondk}) with $\CF'$ in place of $\CF$ holds for all subspaces $0\not=V''\subset V'\subset U_{\rm max}$. We claim that $\CF := \CF'\cup\{V\}$ does the job.

\medskip
Indeed, take any $\BF_q$-subspaces $0\not=V''\subset V'\subset V$. In the case that $V'\subset U_{\rm max}$, condition (\ref{BV:OFcondk}) follows from the induction hypothesis. In the case that $V''\subset U_{\rm max}$ but $V'\not\subset U_{\rm max}$, by the construction of $U_{\rm max}$ we have on the one hand $V''_k \subset U_{{\rm max},k} \subset E_V$, and on the other hand $V'_k\not\subset E_V$ and hence $V_k = E_V + V'_k$. The last equation implies that the homomorphism of $1$-dimensional $k$-vector spaces $V'_k/E_{V'} \rightarrow V_k/E_V$ is surjective, hence also injective, and thus $E_{V'} = E_V \cap V'_k$.
Together with $V''_k \subset U_{{\rm max},k} \subset E_V$ this implies that $V''_k \subset E_V\cap V'_k = E_{V'}$, which is the first case of (\ref{BV:OFcondk}) with $W=U_{\rm max}$. Finally, in the case that $V''\not\subset U_{\rm max}$, we have $V''_k\not\subset E_V$, thus $V''_k\not\subset E_{V'}$, and hence $V'_k = E_{V'} + V''_k$. This is the second case of (\ref{BV:OFcondk}), where indeed no $W$ with the indicated properties exists, because the only $W\in\CF$ with $V''\subset W$ is $W=V$. Thus $\CF$ has the desired properties.
\end{Proof}

\begin{Cor}
\label{BV:geometry}
\begin{itemize}
\item[(a)]  As a set $B_\CF$ is the union of $\Omega_{\CF'}$ for all flags $\CF'$ of~$V$ with $\CF \subset \CF'$.
\item[(b)]  As a set $U_\CF$ is the union of $\Omega_{\CF'}$ for all flags $\CF'$ of~$V$ with $\CF' \subset \CF$.
\item[(c)] The $U_\CF$ for all flags $\CF$ of~$V$ form an open covering of~$B_V$. 
\end{itemize}
\end{Cor}

\begin{Proof}
Combine Lemma \ref{BV:flagcomp} with Theorem \ref{BV:stratification}.
\end{Proof}


\medskip
Next we want to simplify the description of the strata neighborhoods~$U_\CF$. Write $\CF=\{V_0,\ldots,V_m\}$ for subspaces $0=V_0 \subsetneqq V_1\subsetneqq\ldots\subsetneqq V_m = V$.

\begin{Lem}
\label{BV:UFlemma1}
The forgetful map $\CE_\bullet \mapsto (\CE_{V_1},\ldots,\CE_{V_m})$ induces an isomorphism from $U_\CF$ to the locally closed subscheme $U^\flat_\CF \subset \prod_{i=1}^m P_{V_i}$ representing tuples satisfying:
\begin{itemize}
\item[(a)] \ $\CE_{V_1} \subset \ldots \subset \CE_{V_m}$ and
\item[(b)] \ $V_i \otimes \CO_S = \CE_{V_i} + (V' \otimes \CO_S)$ \ \ \ \ \ \ for all $i$ and all $0 \neq V' \subset V_i$ with $V' \not\subset V_{i-1}$.
\end{itemize}
\end{Lem}

\begin{Proof}
The conditions (a) and (b) are special cases of (\ref{BV:BVcond}) and (\ref{BV:UFcond}); hence the map induces a morphism $U_\CF\to U^\flat_\CF$. To construct a morphism in the other direction consider any tuple $(\CE_{V_1},\ldots,\CE_{V_m})$ over a scheme $S$ that satisfies (a) and (b). For any non-zero subspace $V'\subset\nobreak V$ let $i$ be the unique integer such that $V' \subset V_i$ and $V' \not\subset V_{i-1}$, and set $\CD_{V'} := \CE_{V_i} \cap (V' \otimes \CO_S)$. Then condition (b) implies that 
$$(V' \otimes \CO_S)/\CD_{V'} \ \cong\ \bigl(\CE_{V_i}+(V' \otimes \CO_S)\bigr)/\CE_{V_i} 
\ \cong\ (V_i \otimes \CO_S)/\CE_{V_i}.$$
By the assumption on $\CE_{V_i}$ the right hand side is locally free of rank~$1$; hence so is the left hand side, and so $\CD_{V'}$ defines an $S$-valued point of~$P_{V'}$. Next take another subspace $0\not=V''\subset V'$ and let $j$ be such that $V'' \subset V_j$ and $V'' \not\subset V_{j-1}$. Then $j\le i$ and hence 
$$\CD_{V''}\ :=\ \CE_{V_j} \cap (V'' \otimes \CO_S)\ \subset\ \CE_{V_i} \cap (V' \otimes \CO_S)\ =:\ \CD_{V'}.$$
Thus the tuple $(\CD_{V'})_{V'}$ satisfies the condition (\ref{BV:BVcond}) and defines an $S$-valued point of~$B_V$. Furthermore assume that there exists no $W \in \CF$ with $V'' \subset W$ and $V' \not\subset W$. Then this assumption holds in particular for $W=V_{i-1}$, and since $V' \not\subset V_{i-1}$ we deduce that $V'' \not\subset V_{i-1}$. This implies that $j\ge i$ and hence $j=i$. By condition (b) we therefore have $V_i \otimes \CO_S = \CE_{V_i} + (V'' \otimes \CO_S)$. Intersecting this equation with $V' \otimes \CO_S$ yields
\begin{eqnarray*}
V' \otimes \CO_S 
&=& \bigl(\CE_{V_i} + (V'' \otimes \CO_S) \bigr) \cap (V' \otimes \CO_S) \\
&=& \bigl(\CE_{V_i} \cap (V' \otimes \CO_S) \bigr) + (V'' \otimes \CO_S) \\ 
&=& \CD_{V'} + (V'' \otimes \CO_S).
\end{eqnarray*}
This means that the tuple $(\CD_{V'})_{V'}$ satisfies the condition (\ref{BV:UFcond}) and defines an $S$-valued point of~$U_\CF$. Altogether the construction yields a morphism $U^\flat_\CF\to U_\CF$.
 
\medskip
The construction immediately shows that $\CD_{V_i}=\CE_{V_i}$ for all $1\le i\le m$; hence the composite of $U^\flat_\CF\to U_\CF\to U^\flat_\CF$ is the identity. To show that the composite of $U_\CF\to U^\flat_\CF\to U_\CF$ is the identity consider any tuple $\CE_\bullet \in U_{\CF}(S)$. We must verify that $\CE_{V'} = \CD_{V'} := \CE_{V_i} \cap (V' \otimes \CO_S)$ for any $0 \neq V' \subset V_i$ with $V' \not\subset V_{i-1}$. But $\CE_{V'} \subset \CE_{V_i}$ implies that $\CE_{V'} \subset \CD_{V'}$, which yields a natural surjection $(V' \otimes \CO_S)/\CE_{V'} \onto (V' \otimes \CO_S)/\CD_{V'}$. As both sheaves are locally free of rank~$1$ this surjection is in fact an isomorphism, whence $\CE_{V'} = \CD_{V'}$, as desired. Thus the morphism $U_\CF\to U^\flat_\CF$ has a two-sided inverse and is therefore an isomorphism.
\end{Proof}


\begin{Prop}
\label{BV:OVprop}
For the trivial flag $\CF_0:=\{0,V\}$ we have $\Omega_{\CF_0} = U_{\CF_0}$ and a natural isomorphism 
$$\Omega_{\CF_0} \ \stackrel{\sim}{\longto}\ \Omega_V,\ \ \CE_\bullet \mapsto\CE_V.$$
We identify $\Omega_{\CF_0}$ with $\Omega_V$ through this isomorphism.
\end{Prop}

\begin{Proof}
The first assertion is a special case of Corollary \ref{BV:geometry} (b). Thus Lemma \ref{BV:UFlemma1} for $\CF=\CF_0$ yields an isomorphism from $\Omega_{\CF_0}$ to the open subscheme $U^\flat_{\CF_0}\subset P_V$ representing all $\CE_V$ satisfying $V \otimes \CO_S = \CE_V + (V' \otimes \CO_S)$ for all $0 \neq V' \subset V$. By Proposition \ref{BV:PV2} this subscheme is just~$\Omega_V$.
\end{Proof}


\begin{Lem}
\label{BV:UFlemma2}
The isomorphism in Lemma \ref{BV:UFlemma1} identifies the open subscheme $\Omega_V \subset U_{\CF}$ with the subscheme $\Omega^\flat_V\subset U^\flat_\CF$ representing tuples $(\CE_{V_1},\ldots,\CE_{V_m})$ which in addition satisfy:
\begin{itemize}
\item[(c)] \ $V_i \otimes \CO_S = \CE_{V_i} + (V_{i-1} \otimes \CO_S)$ \ \ \ \ for all\ \ $2\le i\le m$.
\end{itemize}
\end{Lem}


\begin{Proof}
First one easily shows by induction on $i-j$ that (c) is equivalent to 
\begin{itemize}
\item[(c$'$)] \ $V_i \otimes \CO_S = \CE_{V_i} + (V_j \otimes \CO_S)$ \ \ \ \ for all $1\le j\le i \le m$.
\end{itemize}
Next it is immediate from (\ref{BV:UFcond}) and Proposition \ref{BV:OVprop} that the image of $\Omega_V=\Omega_{\CF_0}=U_{\CF_0}$ satisfies (c$'$). To obtain the desired isomorphism it therefore suffices to show that any tuple $\CE_\bullet \in U_{\CF}(S)$ satisfying (c$'$) already lies in $U_{\CF_0}(S)$. For this we must prove that $V' \otimes \CO_S = \CE_{V'} + (V'' \otimes \CO_S)$ for arbitrary $0 \neq V'' \subset V' \subset V$. 

\medskip
Let $i$ be the integer such that $V' \subset V_i$ and $V' \not\subset V_{i-1}$, and let $j$ be the integer such that $V'' \subset V_j$ and $V'' \not\subset V_{j-1}$. Then $j\le i$, and so by (c$'$) and \ref{BV:UFlemma2} (b) and (a) we deduce that
\begin{eqnarray*}
V_i \otimes \CO_S
&=& \CE_{V_i} + (V_{j} \otimes \CO_S) \\
&=& \CE_{V_i} + \CE_{V_j} + (V'' \otimes \CO_S) \\
&=& \CE_{V_i} + (V'' \otimes \CO_S).
\end{eqnarray*}
Intersecting this with $V' \otimes \CO_S$ yields
\begin{eqnarray*}
V' \otimes \CO_S
&=& \bigl( \CE_{V_i} + (V'' \otimes \CO_S) \bigr) \cap (V' \otimes \CO_S) \\
&=& \bigl( \CE_{V_i} \cap (V' \otimes \CO_S) \bigr) + (V'' \otimes \CO_S).
\end{eqnarray*}
But in the proof of Lemma \ref{BV:UFlemma2} we showed that $\CE_{V'} = \CE_{V_i} \cap (V' \otimes \CO_S)$, and so the right hand side is $\CE_{V'} + (V'' \otimes \CO_S)$, as desired.
\end{Proof}


\begin{Prop}
\label{BV:UFprop}
Let $r:=\dim V$. There exists an open embedding $U_\CF \into \BA^{r-1}_{\BF_q}$ such that the boundary $U_\CF\setminus\Omega_V$ is the inverse image of the union of all coordinate hyperplanes.
\end{Prop}

\begin{Proof}
By Lemma \ref{BV:flagcomp} (a) it suffices to prove this when $\CF$ is a complete flag, i.e., when $\dim V_i=i$ for all $1\le i\le m$ and $m=r$. We can then choose a basis $X_1, \ldots, X_r$ of~$V$ such that each $V_i$ is generated by $X_1, \ldots, X_i$.
We abbreviate $W_i := \BF_q X_i$.

\medskip
By Lemmas \ref{BV:UFlemma1} and \ref{BV:UFlemma2} it suffices to prove the assertion for $\Omega^\flat_V \subset U^\flat_\CF$ in place of $\Omega_V \subset U_\CF$. Consider the locally closed subscheme $U^\sharp_\CF \subset \prod_{i=1}^m P_{V_i}$ representing tuples $(\CE_{V_1},\ldots,\CE_{V_m})$ that satisfy
\begin{itemize}
\item[(a)] \ $\CE_{V_1} \subset \ldots \subset \CE_{V_m}$ and
\item[(b$'$)] \ $V_i \otimes \CO_S = \CE_{V_i} + (W_i \otimes \CO_S)$ \ \ \ \ \ \ for all $1\le i\le r$.
\end{itemize}
Here (a) coincides with \ref{BV:UFlemma1} (a), and (b$'$) consists of special cases of the open condition \ref{BV:UFlemma1}~(b); hence $U_\CF^\flat$ is an open subscheme of~$U_\CF^\sharp$. Also, let $\Omega^\sharp_V$ denote the open subscheme of $U_\CF^\sharp$ determined by the condition 
\begin{itemize}
\item[(c)] \ $V_i \otimes \CO_S = \CE_{V_i} + (V_{i-1} \otimes \CO_S)$ \ \ \ \ for all\ \ $2\le i\le m$.
\end{itemize}
Then Lemma \ref{BV:UFlemma2} shows that $\Omega^\flat_V = \Omega^\sharp_V \cap U^\flat_V$. Thus it suffices to prove the assertion for $\Omega^\sharp_V \subset U^\sharp_\CF$ in place of $\Omega^\flat_V \subset U^\flat_\CF$.
We will achieve this by producing an \emph{isomorphism} $U^\sharp_\CF \cong \BA^{r-1}_{\BF_q}$ under which $U^\sharp_\CF\setminus\Omega^\sharp_V$ corresponds to the union of all coordinate hyperplanes. For this note that $S$-valued points of $\BA^{r-1}_{\BF_q}$ amount to $(r-1)$-tuples of sections in $\Gamma(S,\CO_S)$. The desired isomorphism thus results from the following lemma:

\begin{Lem}
\label{BV:UFproplem}
\begin{itemize}
\item[(i)] Consider any sections $a_1,\ldots,a_{r-1}\in\Gamma(S,\CO_S)$. For all $1\le i\le r$ let $\CE_{V_i}$ be the locally free coherent subsheaf of $V_i\otimes\CO_S$ generated by 
$$\bigl\{ X_j\otimes1+X_{j+1}\otimes a_j \bigm| 1\le j<i\bigr\}.$$ 
Then the tuple $(\CE_{V_1},\ldots,\CE_{V_m})$ defines an $S$-valued point of $U^\sharp_\CF$.
\item[(ii)] Every $S$-valued point of $U^\sharp_\CF$ arises as in (i) from unique sections $a_1,\ldots,a_{r-1}$.
\item[(iii)] The tuple defines an $S$-valued point of $\Omega^\sharp_V$ if and only if $a_1,\ldots,a_{r-1}\in\Gamma(S,\CO_S^\times)$.
\end{itemize}
\end{Lem}

To prove this consider first the situation of (i). Then condition (a) is obvious. Also, for any $1\le i\le r$ the set
$$\bigl\{ X_j\otimes1+X_{j+1}\otimes a_j \bigm| 1\le j<i\bigr\} \,\cup\, \bigl\{X_i\otimes1\bigr\}$$
is a basis of $V_i\otimes\CO_S$, because it can be obtained by applying a unipotent matrix to the standard basis $\{X_j\otimes1\,|\,1\,{\le}\, j\,{\le}\,i\}$. Thus the definition of $\CE_{V_i}$ in (i) implies that
$V_i \otimes \CO_S = \CE_{V_i} \oplus (W_i \otimes \CO_S)$. In particular this shows (b$'$) and that $(V_i \otimes \CO_S)/\CE_{V_i} \cong W_i \otimes \CO_S$ is (locally) free of rank~$1$, proving~(i).

\medskip
To prove (ii) let $(\CE_{V_1},\ldots,\CE_{V_m})$ be any $S$-valued point of~$U^\sharp_\CF$. Then for each $1\le i\le r$ the identity induces a surjective homomorphism $W_i \otimes \CO_S \to (V_i \otimes \CO_S)/\CE_{V_i}$ by (b$'$). As both sides are locally free of rank~$1$, this homomorphism is in fact an isomorphism, and so (b$'$) can be strengthened to $V_i \otimes \CO_S = \CE_{V_i} \oplus (W_i \otimes \CO_S)$. Since $V_i=V_{i-1}\oplus W_i$, this shows that $\CE_{V_i}$ is the graph of an $\CO_S$-linear homomorphism $V_{i-1}\otimes\CO_S \to W_i\otimes\CO_S$. In particular we have $\CE_{V_1}=0$. For $2\le i\le r$ the homomorphism sends $X_{i-1}\otimes1$ to $X_i\otimes a_{i-1}$ for a unique section $a_{i-1}\in\Gamma(S,\CO_S)$. In other words there is a unique section $a_{i-1}\in\Gamma(S,\CO_S)$ such that $X_{i-1}\otimes1+X_i\otimes a_{i-1} \in \Gamma(S,\CE_{V_i})$. Since $\CE_{V_{i-1}}\subset \CE_{V_i}$ by~(a), varying $i$ yields unique sections $a_1,\ldots,a_{r-1}\in\Gamma(S,\CO_S)$ such that $X_j\otimes1+X_{j+1}\otimes a_j$ is a section of $\CE_{V_i}$ for all $1\le j<i\le r$. But for fixed~$i$, the proof of (i) shows that these sections for $1\le j<i$ already generate a coherent subsheaf of $V_i\otimes\CO_S$ which is a direct complement of $W_i\otimes\CO_S$. Thus these sections generate $\CE_{V_i}$, proving (ii).

\medskip
Finally, by construction the image of $\CE_{V_i}$ in $(V_i \otimes \CO_S)/(V_{i-1} \otimes \CO_S) \cong W_i\otimes\CO_S$ is the coherent subsheaf generated by $X_i\otimes a_{i-1}$. Thus $a_{i-1}$ is invertible if and only if this image is equal to $W_i\otimes\CO_S$, which in turn is equivalent to the condition~(c). This proves (iii); hence it finishes the proof of Lemma \ref{BV:UFproplem} and of Proposition~\ref{BV:UFprop}.
\end{Proof}


\begin{Thm}
\label{BV:smooth}
The scheme $B_V$ is an irreducible smooth projective variety, and the boundary $B_V \setminus \Omega_V$ is a divisor with normal crossings.
\end{Thm}

\begin{Proof}
Being closed in a projective scheme $B_V$ is projective. By Corollary \ref{BV:geometry} the $U_\CF$ form an open covering, and all of them contain the open stratum~$\Omega_V$. Proposition \ref{BV:UFprop} implies that the $U_\CF$ are irreducible. Together this implies that $B_V$ is irreducible. The remaining assertions also follow from Proposition \ref{BV:UFprop}.
\end{Proof}


\begin{Thm}
\label{BV:maps}
There exist morphisms $\pi_P$ and $\pi_Q$ making the following diagram commute:
$$\xymatrix@+10pt{
P_V   &  B_V \ar[l]_-{\pi_P} \ar[r]^-{\pi_Q}  &  Q_V \\
&   \Omega_V \ar@{^{ (}->}[ul]^{\ref{mod:omega1}} 
             \ar@{^{ (}->}[ur]_{\ref{mod:omega2}} 
             \ar@{^{ (}->}[u]^{\ref{BV:OVprop}}   &   \\}$$
\end{Thm}

\begin{Proof}
The map $\CE_\bullet \mapsto \CE_{V}$ clearly induces a morphism $\pi_P: B_V \to P_V$ which by Proposition \ref{BV:OVprop} makes the triangle on the left hand side commute. To construct $\pi_Q$ we use the modular interpretation of $Q_V$ from Theorem \ref{mod:QV} and associate to any $S$-valued point $\CE_\bullet$ of $B_V$ an $S$-valued point $(\CL, \rho)$ of $Q_V$ as follows. 

\medskip
Consider the commutative diagram of invertible sheaves $(V' \otimes \CO_S) / \CE_{V'}$ for all $0 \neq V' \subset \nobreak V$ with the natural homomorphisms $(V'' \otimes \CO_S) / \CE_{V''} \to (V' \otimes \CO_S) / \CE_{V'}$ induced by the inclusions $0 \neq V'' \subset V' \subset V$. Dualizing yields a commutative diagram of invertible sheaves ${((V' \otimes \CO_S )}/ \CE_{V'})^{-1}$ with homomorphisms 
$$\psi_{V'}^{V''}: \bigl((V' \otimes \CO_S) / \CE_{V'} \bigr)^{-1} \longrightarrow 
\bigl((V'' \otimes \CO_S) / \CE_{V''}\bigr)^{-1}.$$
We view this diagram as a direct system, but note that it is not filtered, for example because it contains no arrows out of the objects with $\dim V''=1$.
We define $\CL$ as the direct limit (i.e., colimit) of this system.
For any $v \in\!\circV\,$ let $\ell_v: \BF_q v \otimes \CO_S = (\BF_q v \otimes \CO_S)/\CE_{\BF_qv} \to \CO_S$ be the $\CO_S$-linear homomorphism $v \otimes a\mapsto a$. Then $\ell_v$ is a global section of $((\BF_q v \otimes \CO_S)/\CE_{\BF_qv})^{-1}$, and we define $\rho (v)$ as the image of $\ell_v$ in $\Gamma(S, \CL)$ under the natural map from $((\BF_q v \otimes \CO_S)/\CE_{\BF_qv})^{-1}$ to the direct limit~$\CL$. 

We will prove that $\CL$ is an invertible sheaf on $S$ and that $\rho: \circV \to \Gamma(S, \CL)$ is a fiberwise non-zero reciprocal map. Both assertions are Zariski local on~$S$; hence by Corollary \ref{BV:geometry} we may assume that $\CE_\bullet \in U_{\CF}(S)$ for some flag $\CF$ of~$V$, not necessarily maximal. As before we write $\CF=\{V_0,\ldots,V_m\}$ for subspaces $0=V_0 \subsetneqq V_1\subsetneqq\ldots\subsetneqq V_m = V$.

\medskip
We claim that in this case $\CL = ((V_1 \otimes \CO_S)/\CE_{V_1})^{-1}$. To see this, for any $0\not=V'\subset V$ let $i$ be the unique integer such that $V' \subset V_i$ and $V' \not\subset V_{i-1}$. Then (\ref{BV:UFcond}) for the inclusion $V'\subset V_i$ implies that $\psi_{V_i}^{V'}$ is an isomorphism. Moreover, for any $0\not=V''\subset V'$ and $V'' \subset V_j\subset V_i$ with $V'' \not\subset V_{j-1}$ we have $\psi^{V''}_{V'}\circ\psi^{V'}_{V_i} = \psi^{V''}_{V_j}\circ\psi^{V_j}_{V_i}$ where $\psi_{V_j}^{V''}$ is again an isomorphism. This allows us to eliminate all objects except those associated to $V_1,\ldots,V_m$ from the diagram, without changing the direct limit. Afterwards the system is filtered with the final object $((V_1 \otimes \CO_S)/\CE_{V_1})^{-1}$, which is therefore the direct limit, as claimed.

\medskip
In particular the claim implies that $\CL$ is an invertible sheaf. Next, for any $v\in\circV_1$ the isomorphism $\ell_v: \BF_qv \otimes \CO_S \stackrel{\sim}{\to} \CO_S$ corresponds to a nowhere vanishing section of ${(\BF_qv \otimes \CO_S)^{-1}} \allowbreak = {((V_1 \otimes \CO_S)/\CE_{V_1})^{-1}}$.
Thus the claim implies that $\rho(v)$ is a nowhere vanishing section of~$\CL$, and so $\rho$ is fiberwise non-zero. 

\medskip
A direct proof that $\rho$ is reciprocal would be awkward in this general setting. Instead observe that by pullback it suffices to prove this for the universal family over~$U_\CF$. Since $U_\CF$ is reduced by Proposition \ref{BV:UFprop}, it then suffices to prove the identities in \ref{mod:reciprocal} over the dense open subscheme $U_{\CF_0}\cong\Omega_V$. In other words we can now assume that $\CF=\CF_0=\{0,V\}$. Then $V_1=V$ and all $\psi_{V}^{V'}$ are isomorphisms. Thus for each $v\in\circV$ the section $\rho(v) \in \Gamma(S,\CL)$ is by definition the image of the section $1$ under the isomorphisms
$$\vcenter{\vskip-5pt\xymatrix@C+20pt{
\CO_S \ar[r]^-{1\mapsto\ell_v}_-\sim
& (\BF_qv\otimes\CO_S)^{-1} \ar@{=}[r]^-{\psi_{V}^{\BF_qv}}_-\sim
& \bigl((V\otimes\CO_S)/\CE_V\bigr)^{-1} \rlap{\mbox{$\ =\ \CL.$}} \\}}$$
In particular it vanishes nowhere, and its reciprocal $\rho(v)^{-1} \in \Gamma(S,\CL^{-1})$ is the image of $1$ under the isomorphisms
$$\vcenter{\vskip-5pt\xymatrix@C+20pt@R-24pt{
\CO_S \ar@{<-}[r]^-{\ell_v}_-\sim
& \BF_qv\otimes\CO_S \ar@{=}[r]_-\sim
& (V\otimes\CO_S)/\CE_V \rlap{\mbox{$\ =\ \CL^{-1},$}} \\
\ a\ & \ v\otimes a\ \ar@{|->}[l] \ar@{|->}[r] & \ [v\otimes a]\rlap{.}\ }}$$
Thus if $\lambda$ denotes the $\BF_q$-linear map $V\to\Gamma(S,\CL^{-1})$, $v\mapsto [v\otimes1]$, it follows that $\rho$ is the reciprocal of $\lambda$ according to Proposition \ref{mod:rec-lin-1} and hence a reciprocal map, as desired.

\medskip
To summarize we have associated to any $S$-valued point $\CE_\bullet$ of $B_V$ an $S$-valued point $(\CL,\rho)$ of~$Q_V$. As this construction commutes with pullback, it defines a morphism $\pi_Q: {B_V\to Q_V}$. 
Also, over $\Omega_V\subset B_V$ we have seen that $(\CL,\rho)$ is the reciprocal of $(\CL^{-1},\lambda)$. But $(\CL^{-1},\lambda)$ is just the pair corresponding to $\CE_V = \pi_P(\CE_\bullet)$ under the equivalent modular interpretations \ref{mod:PV} and \ref{BV:PV1} of~$P_V$. Thus $\pi_Q|\Omega_V$ is simply the original embedding $\Omega_V\into Q_V$, showing that the triangle on the right hand side of the diagram commutes.
\end{Proof}


\begin{Prop}
\label{BV:BFprop}
There is a natural isomorphism
$$B_\CF \ \stackrel{\sim}{\longto}\  B_{V_1/V_0} \times \ldots \times B_{V_m/V_{m-1}},\ \ 
\CE_\bullet \mapsto \Bigl( \bigl( \CE_{V'}/(V_{i-1}\otimes\CO_S) \bigr)_{V'/V_{i-1}} \Bigr)_{i=1}^m,$$
where $V'$ runs through all subspaces $V_{i-1}\subsetneqq V'\subset V_i$. 
\end{Prop}

\begin{Proof}
Denote the map of the Proposition by $\mu$. We first show that $\mu$ is well-defined. By taking $V'' := V_{i-1}$ and $W := V_{i-1}$ in condition (\ref{BV:BFcond}) we see that indeed $V_{i-1} \otimes \CO_S \subset \CE_{V'}$ for any $V_{i-1}\subsetneqq V'\subset V_i$. Furthermore, every tuple $(\CE_{V'}/(V_{i-1}\otimes\CO_S))_{V'/V_{i-1}}$ clearly satisfies condition (\ref{BV:BVcond}) applied to $B_{V_i/V_{i-1}}$. Finally, the quotient of $(V'/ V_{i-1}) \otimes \CO_S$ by $\CE_{V'}/(V_{i-1}\otimes\CO_S)$ is isomorphic to $(V' \otimes\CO_S)/\CE_{V'}$ and thus again locally free of rank $1$. Hence $\mu$ is well-defined. Clearly $\mu$ is functorial in $\CE_\bullet$; hence it defines a morphism of schemes. 

\medskip
Next we construct a morphism in the opposite direction $\nu: B_{V_1/V_0} \times \ldots \times B_{V_m/V_{m-1}} \rightarrow B_{\CF}$. To a collection of tuples $\CE^i_\bullet$ in $B_{V_i/V_{i-1}}$ for $1\le i\le m$ we assign a tuple $\CD_\bullet$
in $B_\CF$ as follows. For any non-zero subspace $V' \subset V$, let $i$ be the unique integer such that $V' \subset V_i$ and $V' \not\subset V_{i-1}$, and let $\pi: V' \otimes \CO_S \onto (V' + V_{i-1}/V_{i-1}) \otimes \CO_S$ denote the natural surjection. Then we set $\CD_{V'} := \pi^{-1}(\CE^i_{V'+V_{i-1}/V_{i-1}})$. We now verify that the collection $\CD_\bullet$ obtained in this way indeed defines an $S$-valued point of $B_\CF$.

\medskip
First note that because $\pi$ is surjective, the quotient of $V' \otimes \CO_S$ by $\CD_{V'}$ is isomorphic to the quotient of $(V' + V_{i-1}/V_{i-1}) \otimes \CO_S$ by $\CE^i_{V'+V_{i-1}/V_{i-1}}$ and therefore indeed locally free of rank~$1$. Next we prove that $\CD_{V''} \subset \CD_{V'}$ for any non-zero subspace $V'' \subset V'$. By the definition of $\CD_{V'}$ we need to show that $\pi(\CD_{V''}) \subset \CE^i_{V'+V_{i-1}/V_{i-1}}$. If $V'' \subset V_{i-1}$ we have $\pi(V'' \otimes \CO_S) = 0$ and there is nothing to prove. Thus we can assume that $V'' \not\subset V_{i-1}$. But as $V'' \subset V' \subset V_i$, the construction of $\CD_\bullet$ shows that $\CD_{V''}$ is the inverse image of $\CE^i_{V''+V_{i-1}/V_{i-1}}$ under the restriction of $\pi$ to $V'' \otimes \CO_S$. Since $\CE^i_{V''+V_{i-1}/V_{i-1}} \subset \CE^i_{V'+V_{i-1}/V_{i-1}}$ by condition \ref{BV:BVcond} applied to $B_{V_i/V_{i-1}}$, it follows that $\CD_{V''}\subset \CD_{V'}$, as desired. 

\medskip
We have now shown that $\CD_\bullet$ defines an $S$-valued point of $B_V$. To see that this point lies in $B_\CF$ we let $V'$ and $i$ and $\pi$ be as before and let $V'' \subset V'$ be a non-zero subspace such that there exists $W \in \CF$ with $V'' \subset W$ and $V' \not\subset W$. Then from $V' \not\subset W$ and $V' \subset V_i$ we conclude that $W \subset V_{i-1}$. This implies that $\pi(V'' \otimes \CO_S) \subset \pi(W \otimes \CO_S) \subset \pi(V_{i-1} \otimes \CO_S) = 0$ and hence $V'' \otimes \CO_S \subset \CD_{V'}$. By (\ref{BV:BFcond}) this means that $\CD_\bullet$ lies in $B_\CF$, as desired. This finishes the construction of~$\nu$.

\medskip
The definitions of $\mu$ and $\nu$ directly imply that $\mu \circ \nu = \id$. 
To show that $\nu \circ \mu = \id$, let $\CE_\bullet$ be an $S$-valued point of $B_\CF$ and let $\CD_\bullet$ denote its image under $\nu \circ \mu$. To prove that $\CE_\bullet = \CD_\bullet$ we let $0 \neq V' \subset V$ be any non-zero subspace and let $i$ and $\pi$ be as before. Chasing through the constructions of $\mu$ and $\nu$ yields $\CD_{V'} = \pi^{-1}(\CE_{V'+V_{i-1}}/V_{i-1} \otimes \CO_S) = (\CE_{V'+V_{i-1}}) \cap (V' \otimes \CO_S)$. From this and from condition (\ref{BV:BVcond}) we see that $\CE_{V'} \subset \CD_{V'}$, and thus there exists a natural surjection $V' \otimes \CO_S/ \CE_{V'} \onto V' \otimes \CO_S/ \CD_{V'}$. But since both quotient sheaves are locally free of rank $1$, this surjection is an isomorphism, and hence $\CE_{V'} = \CD_{V'}$, as desired.
\end{Proof}


\begin{Prop}
\label{BV:OFprop}
The isomorphism of Proposition \ref{BV:BFprop} induces an isomorphism
$$\Omega_\CF \ \stackrel{\sim}{\longto}\  \Omega_{V_1/V_0} \times \ldots \times \Omega_{V_m/V_{m-1}},\ \ 
\CE_\bullet \mapsto \bigl(\CE_{V_i}/(V_{i-1}\otimes\CO_S)\bigr)_{i=1}^m.$$
\end{Prop}

\begin{Proof}
Throughout the proof we identify $\Omega_{V_i/V_{i-1}}$ with the open stratum $\Omega_{\{0, V_i/V_{i-1}\}}$ of $B_{V_i/V_{i-1}}$, as explained in Proposition \ref{BV:OVprop}. We first prove that the image of $\Omega_\CF$ under $\mu$ is contained in $\Omega_{V_1/V_0} \times \ldots \times \Omega_{V_m/V_{m-1}}$. Let $\CE_\bullet$ be any $S$-valued point of $\Omega_\CF$ and let $V'' \subset V' \subset V$ be non-zero subspaces with $V_{i-1} \subsetneqq V'' \subset V' \subset V_i$. Then we need to show that $(V'/V_{i-1}) \otimes \CO_S = \CE_{V'}/(V_{i-1} \otimes \CO_S) + (V''/V_{i-1}) \otimes \CO_S$. But since there exists no $W \in \CF$ such that $V'' \subset W$ and $V' \not\subset W$ we have $V' \otimes \CO_S = \CE_{V'} + (V'' \otimes \CO_S)$ by condition (\ref{BV:UFcond}). Thus dividing the last equation by $V_{i-1} \otimes \CO_S$ yields the desired result.

\medskip
Next we show that the image of $\Omega_{V_1/V_0} \times \ldots \times \Omega_{V_m/V_{m-1}}$ under $\nu$ is contained in $\Omega_\CF$. Let $(\CE^1_\bullet, \ldots, \CE^m_\bullet)$ be an $S$-valued point of $\Omega_{V_1/V_0} \times \ldots \times \Omega_{V_m/V_{m-1}}$ and let $\CD_\bullet$ denote its image under~$\nu$. Given non-zero subspaces $V'' \subset V' \subset V$ such that there exists no $W \in \CF$ with $V'' \subset W$ and $V' \not\subset W$, we have to show that $V' \otimes \CO_S = \CD_{V'} + (V'' \otimes \CO_S)$. Let~$i$ be the unique integer such that $V' \subset V_i$ and $V' \not\subset V_{i-1}$. Then we automatically have $V'' \subset V_i$ and $V'' \not\subset V_{i-1}$ because otherwise $W:=V_{i-1} \in \CF$ would contain $V''$ but not $V'$. As earlier we denote by $\pi$ the natural surjection $V' \otimes \CO_S \onto (V'+V_{i-1}/V_{i-1}) \otimes \CO_S$. Thus $\CD_{V'} = \pi^{-1} \bigl( \CE^i_{V'+V_{i-1}/V_{i-1}} \bigr)$ by the definition of~$\nu$, and since $\CE^i_\bullet$ satisfies condition (\ref{BV:UFcond}) for the trivial flag of $V_i/V_{i-1}$ we have $(V'+V_{i-1}/V_{i-1}) \otimes \CO_S = \CE^i_{V'+V_{i-1}/V_{i-1}} + (V''+V_{i-1}/V_{i-1}) \otimes \CO_S$. Using this and the fact that $V_{i-1} \otimes \CO_S$ is contained in $\CD_{V'}$ by condition (\ref{BV:BFcond}) we conclude
\begin{eqnarray*}
V' \otimes \CO_S
&=& \pi^{-1} \bigl( (V'+V_{i-1}/V_{i-1}) \otimes \CO_S \bigr)  \\
&=& \pi^{-1} \bigl( \CE^i_{V'+V_{i-1}/V_{i-1}} + (V''+V_{i-1}/V_{i-1}) \otimes \CO_S \bigr) \\ 
&=& \pi^{-1} \bigl( \CE^i_{V'+V_{i-1}/V_{i-1}} \bigr) + \pi^{-1} \bigl( (V''+V_{i-1}/V_{i-1}) \otimes \CO_S \bigr) \\
&=& \CD_{V'} + (V''+V_{i-1}) \otimes \CO_S \\
&=& \CD_{V'} + (V'' \otimes \CO_S),
\end{eqnarray*}
finishing the proof.
\end{Proof}


\begin{Cor}
\label{BV:strataclosure}
The closure of $\Omega_\CF$ in $B_V$ is~$B_\CF$.
\end{Cor}

\begin{Proof}
Combination of Propositions \ref{BV:BFprop} and \ref{BV:OFprop} and the fact that, as a consequence of Theorem~\ref{BV:smooth}, each $\Omega_{V_i/V_{i-1}} \subset B_{V_i/V_{i-1}}$ is dense.
\end{Proof}


\begin{Exer}
In the case $\dim V=3$ the variety $B_V$ is the blowup of $P_V$ in all zero-dimensional strata, i.e., the blowup of $\BP^2_{\BF_q}$ in all $\BF_q$-rational points. It is also the blowup of $Q_V$ in all (reduced) zero-dimensional strata.
\end{Exer}



 
\end{document}